\documentclass[VANCOUVER,STIX1COL]{WileyNJD-v2}

\articletype{Research Article}%

\received{26 April 2016}
\revised{6 June 2016}
\accepted{6 June 2016}

\newcommand{\bs}[1]{\ensuremath{\boldsymbol{#1}}}

\newcommand{\statevar}{u}
\newcommand{\bstatevar}{\bs{u}}
\newcommand{\bstatevarupdate}{\bs{\hat{u}}}
\newcommand{\param}{\rho}
\newcommand{\bparam}{\bs{\rho}}
\newcommand{\bparamupdate}{\bs{\hat{\rho}}}
\newcommand{\paramlbound}{\rho_{\ell}}
\newcommand{\bparamlbound}{\bs{\rho}_{\ell}}
\newcommand{\adjoint}{\lambda}
\newcommand{\badjoint}{\bs{\lambda}}
\newcommand{\badjointupdate}{\bs{\hat{\lambda}}}
\newcommand{\lboundmultiplier}{z_{\ell}}
\newcommand{\blboundmultiplier}{\bs{z}_{\ell}}
\newcommand{\blboundmultiplierupdate}{\bs{\hat{z}}_{\ell}}
\newcommand{\dimstate}{n_{\bstatevar}}
\newcommand{\dimparam}{n_{\bparam}}

\newcommand{\spatialvar}{\bs{y}}

\newcommand{\objfun}{f}
\newcommand{\logbarrierfun}{\varphi}
\newcommand{\constfun}{c}
\usepackage{pgfplots}
\usepgfplotslibrary{colorbrewer}

\newcommand{\logLogSlopeTriangle}[5]
{
	
	\pgfplotsextra
	{
		\pgfkeysgetvalue{/pgfplots/xmin}{\xmin}
		\pgfkeysgetvalue{/pgfplots/xmax}{\xmax}
		\pgfkeysgetvalue{/pgfplots/ymin}{\ymin}
		\pgfkeysgetvalue{/pgfplots/ymax}{\ymax}
		
		\pgfmathsetmacro{\xArel}{#1}
		\pgfmathsetmacro{\yArel}{#3}
		\pgfmathsetmacro{\xBrel}{#1-#2}
		\pgfmathsetmacro{\yBrel}{\yArel}
		\pgfmathsetmacro{\xCrel}{\xArel}
		
		\pgfmathsetmacro{\lnxB}{\xmin*(1-(#1-#2))+\xmax*(#1-#2)} 
		\pgfmathsetmacro{\lnxA}{\xmin*(1-#1)+\xmax*#1} 
		\pgfmathsetmacro{\lnyA}{\ymin*(1-#3)+\ymax*#3} 
		\pgfmathsetmacro{\lnyC}{\lnyA+#4*(\lnxA-\lnxB)}
		\pgfmathsetmacro{\yCrel}{\lnyC-\ymin)/(\ymax-\ymin)}
		\pgfmathsetmacro{\negdydx}{int(-1*(#4))}
		
		\coordinate (A) at (rel axis cs:\xArel,\yArel);
		\coordinate (B) at (rel axis cs:\xBrel,\yBrel);
		\coordinate (C) at (rel axis cs:\xCrel,\yCrel);
		
		\draw[#5]   (A)-- node[pos=0.5,anchor=south] {1}
		(B)-- 
		(C)-- node[pos=0.5,anchor=west] {\negdydx}
		cycle;
	}
}

\newtheorem{innercustomproperty}{Property}
\newenvironment{customproperty}[1]
{\renewcommand\theinnercustomproperty{#1}\innercustomproperty}
{\endinnercustomproperty}

\newtheorem{innercustomremark}{Remark}
\newenvironment{customremark}[1]
{\renewcommand\theinnercustomremark{#1}\innercustomremark}
{\endinnercustomremark}

\newtheorem{innercustomnote}{Note}
\newenvironment{customnote}[1]
{\renewcommand\theinnercustomnote{#1}\innercustomnote}
{\endinnercustomnote}

\begin{document}

\title{A Scalable Interior-Point Gauss-Newton Method for PDE-Constrained Optimization with Bound Constraints}

\author[1]{Tucker Hartland}

\author[1]{Cosmin G. Petra}

\author[2]{No\'{e}mi Petra}

\author[1]{Jingyi Wang}

\authormark{Hartland \textsc{et al}}

\address[1]{\orgdiv{Center for Applied Scientific Computing}, \orgname{Lawrence Livermore National Laboratory}, \orgaddress{\state{California}, \country{USA}}}

\address[2]{\orgdiv{Department of Applied Mathematics}, \orgname{University of California, Merced}, \orgaddress{\state{California}, \country{USA}}}

\corres{*Tucker Hartland, \email{hartland1@llnl.gov}}


\abstract[Abstract]{We present a scalable approach to solve a class of elliptic partial differential equation (PDE)-constrained optimization problems with bound constraints. 
	This approach utilizes a robust full-space interior-point (IP)-Gauss-Newton optimization method. 
	To cope with the poorly-conditioned IP-Gauss-Newton saddle-point linear systems that need to be solved, once per optimization step, we propose two spectrally related preconditioners. These preconditioners leverage the limited informativeness of data in regularized PDE-constrained optimization problems. 
	A block Gauss-Seidel preconditioner is proposed for the GMRES-based solution of the IP-Gauss-Newton linear systems. It is shown, for a large-class of PDE- and bound-constrained optimization problems, that the spectrum of the block Gauss-Seidel preconditioned IP-Gauss-Newton matrix is asymptotically independent of discretization and is not impacted by the ill-conditioning that notoriously plagues interior-point methods.
	We exploit symmetry of the IP-Gauss-Newton linear systems and propose a regularization and log-barrier Hessian preconditioner for the preconditioned conjugate gradient (PCG)-based solution of the related IP-Gauss-Newton-Schur complement linear systems. The eigenvalues of the block Gauss-Seidel preconditioned IP-Gauss-Newton matrix, that are not equal to one, are identical to the eigenvalues of the regularization and log-barrier Hessian preconditioned Schur complement matrix.
	The scalability of the approach is demonstrated on an example problem with bound and nonlinear elliptic PDE  constraints. The numerical solution of the optimization problem is shown to require a discretization independent number of IP-Gauss-Newton linear solves.
	Furthermore, the linear systems are solved in a discretization and IP ill-conditioning independent number of preconditioned Krylov subspace iterations.
	The parallel scalability of preconditioner and linear system matrix applies, achieved with algebraic multigrid based solvers, and the aforementioned algorithmic scalability permits a parallel scalable means to compute solutions of a large class of PDE- and bound-constrained problems.  
}



\maketitle

\footnotetext{\textbf{Abbreviations:} PDE, partial differential equation; IPM, interior-point method; IP, interior-point; GMRES, generalized minimal residual method; MINRES, minimal residual method; CG, conjugate gradient; AMG, algebraic multigrid;}

\section{Introduction}\label{sec:introduction}

Partial differential equation (PDE)-based simulation is an indispensable tool for aiding our understanding of physical processes. However, such models invariably contain uncertainty due to, e.g., unknown parameter fields, source terms, initial conditions, boundary conditions, geometries, modeling errors and errors due to discretization. In this work, we target elliptic PDE- and bound-constrained optimization problems and assume that the PDE-based model uncertainty is due entirely to lack of knowledge of a spatially distributed parameter field. The solution of PDE- and bound-constrained optimization problems is a means to learn, from data, unknown or uncertain aspects of PDE-based models that satisfy additional bound constraints. Bound constraints are natural components of physics-based PDEs.
Examples of application specific bound constraints include, but are not limited to: sign-definiteness of a spatially distributed friction field~\cite{perego2012}, and being constrained to the unit interval in order to represent a material fraction field as in density-based topology optimization~\cite{bendsoe2013}. General inequality constraints can be recast as bound constraints through the introduction of one or more auxiliary slack variables. An example application that includes inequality constraints is the nonnegativity of signed-distances between two or more bodies which is used to describe the mechanical deformation of impenetrable elastic bodies in contact~\cite{wriggers2006}. Bound constraints are, however, the source of additional computational challenges as they introduce nonsmooth complementarity conditions into the Karush-Kuhn-Tucker (KKT)~\cite{nocedalwright2006} necessary conditions for optimality.
An alternative to bound constraints is to re-parametrize the PDE-model, e.g., as in~\cite[Section 2.5]{petra2011}. 
An invertible re-parametrization that eliminates bound constraints must necessarily 
be nonlinear and can lead to additional computational challenges.

In this work, we propose a scalable computational method, appropriate for solving large-scale nonlinear elliptic PDE- and bound-constrained optimization problems on parallel high-performance computing systems. This approach directly and systematically addresses the bound constraints and does not require re-parametrization. For the ``outer'' optimization loop, the proposed method is composed of a robust Newton-based interior-point method with a globalizing filter line-search as detailed in Section~\ref{sec:IPGaussNewton}. Interior-point methods have a number of benefits over active-set and semi-smooth Newton methods~\cite{hintermuller2002primal}, such as avoiding infeasible parameters at intermediate stages of the optimization procedure that can cause various components of a computational framework to breakdown as reported e.g., in~\cite[Chapter 4.2]{hans2017} and~\cite[Section 5.2]{hartland2021ssN}. One disadvantage of Newton-based interior-point methods is that in addition to the challenges of solving large-scale linear systems, the linear systems that arise from an IPM have the undesirable property of
becoming arbitrarily ill-conditioned as the optimizer estimate approaches the optimal point. For this reason, in large-scale settings it is essential to have robust preconditioning strategies.

We utilize a Gauss-Newton variant of a primal-dual interior-point method. 
While the convergence rate of the Gauss-Newton method is less than Newton's method, it does converge rapidly for many problems, such as that studied here. The Gauss-Newton method has the desirable property that it generates search directions that are guaranteed to be of descent when the PDE-constraint is nearly satisfied, which is critical for line-search optimization methods. This is not the case for Newton's method when applied to nonconvex optimization problems. 

Preconditioners are essential for the solution of large-scale IP-Gauss-Newton linear systems by Krylov subspace methods. In Section~\ref{sec:linearalgebra} we discuss two related preconditioning strategies. It is established, for typical PDE-constrained optimization problems, that the eigenvalues of the proposed preconditioned linear system matrices cluster near unity in a manner that is independent of the mesh and the log-barrier parameter. The log-barrier parameter is a notorious source of ill-conditioning in IP-Newton and IP-Gauss-Newton linear systems. The proposed framework is then applied for the numerical solution of a large-scale nonlinear elliptic PDE- and bound-constrained optimization problem as detailed in Section~\ref{sec:problemdescription}. Algorithmic and parallel scaling results are reported in Section~\ref{sec:numericalresults}.

\paragraph{Related work}

Numerous works have discussed computationally efficient methods for large-scale PDE-constrained optimization, such as those wherein Newton~\cite{birosghattas2005}, Gauss-Newton~\cite{akcelik2002, petrazhustadleretal2012} and quasi-Newton~\cite{vuchkov2020, petra2023} based optimization methods are employed. The Newton and Gauss-Newton linear system matrices arise from a linearization of the (perturbed) first-order optimality conditions and are  indefinite saddle-point matrices~\cite{benzi2005numerical}, for which effective preconditioners are challenging to construct and direct methods~\cite{davis2016survey}, which scale poorly, are often applied.
In large-scale settings, effective preconditioners for Krylov subspace solvers of Newton and Gauss-Newton linear systems are essential. 
Preconditioners have been proposed that 
perform well in various regimes, such preconditioners include those~\cite{buighattasetal2013} that exploit the limited-informativeness of the underlying data, and those~\cite{joachim2007, pearson2011, pearson2012, Alger19, alger2024point, cyr2023} that are particularly robust to the parameters that define the regularization component of the objective functional. With the exception of~\cite{pearson2011, pearson2017}, the authors are not aware of a work in which the previously mentioned preconditioners are applied to interior-point Newton-based linear systems that arise from PDE- and bound-constrained optimization problems. While the performance of the preconditioner from~\cite{pearson2011, pearson2017} is robust to the regularization parameters, the nonsingularity requirement of the Hessian of the objective with respect to the state is undesirable. It is undesirable as, for many PDE-constrained optimization problems, the nonsingularity requirement ultimately means that the state must be observed throughout the entire domain, which makes it unsuitable for a significant number of problems. 

Here, the proposed computational framework makes use of a block Gauss-Seidel preconditioner presented in~\cite{pestana2016null},
for IP-Gauss-Newton linear system solution by the GMRES method~\cite{saad1986}.
This preconditioner does not 
require that the Hessian of the objective with respect to the state is nonsingular.
Furthermore, this preconditioner exploits known structure present in a large class of
PDE-constrained optimization problems, wherein the data is insufficiently informative and ultimately makes the associated unregularized problem ill-posed~\cite{vogel2002}. 
We establish, from known properties of this class of PDE-constrained optimization problems~\cite{ghattas2021}, that the
spectrum of the preconditioned linear system matrix does not degrade with respect to mesh refinement and is robust to ill-conditioning due to IPMs.
The linear system that describes the linearized first-order optimality conditions in the reduced-space approach is the Schur complement system of the linearized optimality conditions in the full-space approach, see e.g,~\cite[Chapter 1]{hartland2022thesis}. We exploit this fact to design an intimately related preconditioned CG solver that shares identical spectral properties as the block Gauss-Seidel preconditioned GMRES solver. It is to be noted that, while the focus of this work is line-search based methods, trust-region methods~\cite[Chapter 4]{nocedalwright2006} have been successfully applied to both PDE-constrained optimization problems and PDE- and bound-constrained optimization problems, see e.g.,~\cite{heinkenschloss2014, dennis1998}.

\paragraph{Contributions}

The article makes the following contributions to large-scale PDE- and bound-constrained optimization: (1)~explicating the connection between two related preconditioners and their effectiveness for the IP-Gauss-Newton linear systems in PDE- and bound-constrained optimization; (2)~rigorously showing that for a large class of PDE- and bound-constrained optimization problems that said preconditioners are robust to those elements of IPMs that are notorious sources of ill-conditioning; (3)~demonstrating both algorithmic- and strong-scaling of two Interior-Point-Gauss-Newton-Krylov methods on a nonlinear elliptic PDE- and bound-constrained optimization example problem in a parallel framework using the finite element library MFEM~\cite{andrej2024}.

\section{Preliminaries}\label{sec:preliminaries}

In this section, the target class of regularized PDE- and bound-constrained optimization problems is described. We discuss challenges of reduced- and full-space approaches for solving such problems. In Section~\ref{subsec:interiorpoint} we then overview an interior-point method for discretized PDE- and bound-constrained optimization problems, with a stopping criteria (see Section~\ref{subsec:stoppingcriteria}) that is well defined with respect to an infinite-dimensional formulation of the problem.

\paragraph{Notation}
Bold faced symbols generally indicate non-scalar quantitities, \textit{e.g.}, vector and matrices. $\bs{x}_{i}$ denotes the $i$-th component of the vector $\bs{x}$ and $\bs{A}_{i,j}$ is the entry in the $i$-th row and $j$-th column of the matrix $\bs{A}$. $\mathcal{K}_{m}(\bs{A},\bs{r})$ denotes the Krylov subspace $\text{Span}\lbrace \bs{r},\bs{A}\bs{r},\bs{A}^{2}\bs{r},\dots,\bs{A}^{m-1}\bs{r}\rbrace$, of dimension no greater than $m$. The symbol $\bs{1}$ refers to a vector with all entries equal to one, $\bs{0}$ refers to a vector with all entries equal to zero, and $\bs{I}$ is an identity matrix whose dimension, the authors hope, can be made clear from context. The elementwise Hadamard product is indicated by $\odot$, that is $\left(\bs{x}\odot\bs{y}\right)_{i}=\bs{x}_{i}\bs{y}_{i}$. The standard Euclidean vector norm $\sqrt{\bs{x}^{\top}\bs{x}}$, is denoted by $\|\bs{x}\|_{2}$, whereas weighted inner products, e.g., the $\bs{M}$-weighted inner product $\sqrt{\bs{x}^{\top}\bs{M}\bs{x}}$ is denoted by $\|\bs{x}\|_{\bs{M}}$, for a symmetric and positive definite weight matrix $\bs{M}$. 
$\overline{\mathcal{S}}$ denotes the closure of a set $\mathcal{S}$, and $\emptyset$ is the empty set. The $j$-th largest eigenvalue of a symmetric matrix $\bs{A}$, will be denoted $\lambda_{j}(\bs{A})$. The Hilbert space of functions that are square integrable on a bounded domain $\Omega$ is $L^{2}(\Omega)$, and $H^{1}(\Omega)$ is the Sobolev space that contains all $L^{2}(\Omega)$ functions with weak first derivatives that are also $L^{2}(\Omega)$. The norm associated to the Hilbert space $L^{2}(\Omega)$ is given by $\|u\|_{L^{2}(\Omega)}=\left(\int_{\Omega}(u(\spatialvar))^{2}\mathrm{d}\spatialvar\right)^{1/2}$, for $u\in L^{2}(\Omega)$.

\subsection{Optimization problem formulation}

The primary focus of this work is to determine a spatially distributed state $\statevar(\spatialvar)$, and parameter $\param(\spatialvar)$, by solving a
PDE- and bound-constrained optimization problem of the form
\begin{subequations} \label{eq:infdimproblem}  
	\begin{align}
	\min_{(\statevar,\param)\in\mathcal{V}\times\mathcal{M}}&\objfun(\statevar,\param), \label{eq:objective}\\
	\text{such that } &\constfun(\statevar,\param,\adjoint)=0,\quad\forall \adjoint \in\mathcal{V}_{0}, \label{eq:PDEconstraint}\\
	\text{ and } &\param(\spatialvar) \geq \paramlbound(\spatialvar),\quad \forall \spatialvar \in\overline{\Omega}. \label{eq:boundconstraint}
	\end{align}
\end{subequations}
\noindent Here, $\objfun(\statevar,\param)$ is the objective functional that mathematically encodes 
what is deemed optimal, $\mathcal{V}$ is the space of admissible states, $\mathcal{M}$ is the space of admissible parameters, $\constfun$ is a potentially nonlinear weak form of the partial differential equality constraint, $\paramlbound$ is a pointwise parameter lower-bound, and $\mathcal{V}_{0}$ is the space that is tangent to the affine subspace $\mathcal{V}$, that is $(v+v_{0})\in\mathcal{V}$ for any $v\in\mathcal{V}$ and $v_{0}\in\mathcal{V}_{0}$. Here, $\adjoint\in\mathcal{V}_{0}$ is an arbitrary test function used to describe the PDE, but later,
in the context of the Lagrangian formalism for constrained optimization, it will be
a Lagrange multiplier associated to the PDE-constraint. Throughout this work, it is assumed that the objective functional $\objfun(\statevar,\param)$ separates into data-misfit and regularization components $\objfun(\statevar,\param)=\objfun_{\text{misfit}}(\statevar)+\objfun_{\text{reg}}(\param)$, as is the case for many inverse problems governed by PDEs, see e.g.,~\cite{ghattas2021}. The data-misfit $\objfun_{\text{misfit}}$ is a data-fidelity term that measures discrepancy between data and PDE-model predictions associated to the data. The regularization $\objfun_{\text{reg}}$ penalizes unwanted features of the parameter $\param$, e.g., high frequency oscillations so that the regularization promotes spatially smooth parameter reconstructions.

One means to solve optimization problems with PDE-constraints is the reduced-space approach~\cite{hinze2008, borzi2011}, in which one formally eliminates the PDE-solution $\statevar$, by making use of the implicitly defined parameter $\param$ to PDE-solution $\statevar(\param)$ map, whose existence is guaranteed by uniqueness properties of the PDE. In the reduced-space approach, evaluating various quantities such as the reduced-space functional $f(\param):=f(\statevar(\param),\param)$, or its various derivatives by, e.g., the adjoint method~\cite{gunzburger2002, jacobson1968}, can be computationally expensive as they require solution computations of nonlinear and linearized PDEs. Furthermore, optimizer estimation may require many such evaluations and thus computing the solution of a significant number of PDEs may then be required to solve the optimization problem. However, in the absence of bound constraints the reduced-space approach has a relatively simple optimality system such that auxiliary merit functionals and filters are not needed in order that the associated line-search method is globalized~\cite{nocedalwright2006}. In this work, we use a full-space optimization approach in which objective evaluation and derivative computations are significantly cheaper as they do not require nonlinear and linearized PDE-solution computations. This fact is especially relevant with regard to globalizing inexact backtracking line-search~\cite[Chapter 3]{nocedalwright2006}.

\subsection{IPM and discretized optimality conditions} \label{subsec:interiorpoint}

The proposed approach falls under the umbrella of interior-point methods. Specifically, we employ a filter line-search interior-point method (IPM)~\cite{wachter2006} that has emerged as one of the most robust methods for nonlinear nonconvex optimization. This method also possesses best-in-class global and local convergence properties~\cite{wachter2005, waecther_05_ipopt2}. An IPM involves solving a sequence of log-barrier subproblems
\begin{subequations} \label{eq:infdimlogbarrier} 
	\begin{align}
	\min_{(\statevar,\param)\in\mathcal{V}\times\mathcal{M}}\logbarrierfun(\statevar,\param)&:=\objfun(\statevar,\param)-\mu\int_{\Omega}\log(\param-\paramlbound)\,\mathrm{d}\bs{x}, \label{eq:logbarobjective} \\
	\text{ such that } c(\statevar,\param,\adjoint)&=0,\,\,\,\forall \adjoint \in\mathcal{V}_{0}, \label{eq:logbareqcons}
	\end{align}
\end{subequations}
whose optimality conditions can be viewed as smoothed versions of the nonsmooth optimality conditions associated to Equation~\eqref{eq:infdimproblem}. Mathematically, a strictly positive so-called log-barrier parameter $\mu$ is introduced and a sequence of problems for $\mu\rightarrow 0^+$, are solved inexactly. The spatially distributed inequality constraint $ \param \geq \paramlbound$ is implicitly enforced in the form of $\param > \paramlbound$ via a line-search that respects a fraction-to-boundary rule.

In order to obtain a computational solution of the problem described by Equation~\eqref{eq:infdimlogbarrier}, said equation must be discretized. In this work, the discretization of the PDE-constraint, denoted in Equation~\eqref{eq:logbareqcons}, as well as all continuous fields are obtained by the finite element method, but other discretization methods such as finite differences~\cite{strikwerda2004} can be used as well. It is essential to ensure that the discretized problem, including norms, inner products and stopping criteria (see Section~\ref{subsec:stoppingcriteria}) are consistent with the underlying Hilbert spaces. If these spaces are not properly accounted for then the IPM is prone to mesh dependent performance~\cite{petra2023, vuchkov2020}, such as requiring a mesh dependent number of optimization steps that is ultimately detrimental to the scaling of the method with respect to mesh refinement. This is especially pervasive when the meshes are nonuniform, for example, as the result of adaptive mesh refinement. We refer the reader to~\cite{petra2023} for details and a comparative discussion of the differences between the discretized IPM of this work and the finite element oblivious IPM from~\cite{wachter2006}. 

Upon discretization by finite elements, Equation~\eqref{eq:infdimlogbarrier} becomes 
\begin{subequations} \label{eq:finitedimlogbarrier}
	\begin{align}
	&\min_{(\bstatevar,\bparam)\in\mathbb{R}^{\dimstate}\times\mathbb{R}^{\dimparam}}\logbarrierfun_{h}(\bstatevar,\bparam):=\objfun_{h}(\bstatevar,\bparam)-\mu\bs{1}^{\top}\bs{M}\log(\bparam-\bparamlbound), \label{eq:logbarobjective_discrete}\\
	&\text{ such that }\bs{c}(\bstatevar,\bparam)=\bs{0}.\label{eq:logbareqcons_discrete}
	\end{align}
\end{subequations}
Here, the finite element approximation $(\statevar_{h}(\bs{x})$, $\param_{h}(\bs{x}))\in\mathcal{V}_{h}\times \mathcal{M}_{h}$ of $(\statevar(\bs{x})$, $\param(\bs{x}))\in\mathcal{V}\times \mathcal{M}$ is represented by vectors $\bstatevar\in\mathbb{R}^{\dimstate}, \bparam\in\mathbb{R}^{\dimparam}$ via basis elements 
$\lbrace \phi_{i}(\bs{x})\rbrace_{i=1}^{\dimstate}\subset \mathcal{V}_{h}$, $\lbrace \psi_{i}(\bs{x})\rbrace_{i=1}^{\dimparam}\subset \mathcal{M}_{h}$ as
\begin{align}
\statevar_{h}(\bs{x})=\sum_{i=1}^{\dimstate}\bstatevar_{i}\phi_{i}(\bs{x}),\text{       }
\param_{h}(\bs{x})=\sum_{i=1}^{\dimparam}\bparam_{i}\psi_{i}(\bs{x}),
\end{align}
furthermore $\bs{M}_{i,j}=\int_{\Omega}\psi_{i}(\bs{x})\,\psi_{j}(\bs{x})\,\mathrm{d}\bs{x}$ is the mass matrix with respect to the space $\mathcal{M}_{h}$ that is used to discretize the parameter $\param(\bs{x})$. The discretized partial differential equality constraint $\bs{\constfun}$ is given by
\begin{align*}
\left[\bs{\constfun}(\bstatevar,\bparam)\right]_{i}&=c(\statevar_{h},\param_{h},\phi_{i}),\quad 1\leq i\leq \dimstate.	
\end{align*}
It is to be noted that $\log(\bparam-\bparamlbound)$ is the logarithm of the difference of nodal discretizations of $\param$ and $\paramlbound$, rather than
a vector representation of a nodal discretization of $\log(\param-\paramlbound)$. The two discretized quantities converge to one another in the continuous limit, but we prefer the former due to its simplicity. Formally, the IPM requires the solution of a sequence of subproblems (see Equation~\eqref{eq:finitedimlogbarrier}), each said subproblem solution necessarily satisfies the first-order optimality conditions~\cite{nocedalwright2006}
\begin{subequations} \label{eq:optconditionslogbarrier} 
	\begin{align}
	\bs{\nabla}_{\! \bstatevar}\logbarrierfun_{h}+\bs{J}_{\! \bstatevar}^{\top}\badjoint
	&=\bs{\nabla}_{\! \bstatevar}\objfun_{h}+\bs{J}_{\! \bstatevar}^{\top}\badjoint=\bs{0}, \label{eq:discretestationarityu} \\
	\bs{\nabla}_{\! \bparam}\logbarrierfun_{h}+\bs{J}_{\! \bparam}^{\top}\badjoint&=\bs{\nabla}_{\! \bparam}\objfun_{h}+\bs{J}_{\! \bparam}^{\top}\badjoint-
	\bs{M}_{\! L}(\mu\bs{1}/(\bparam-\bparamlbound))=\bs{0}, \label{eq:discretestationaritym} \\
	\bs{c}&=\bs{0}, \label{eq:discreteequalityconstraint}
	\end{align}
\end{subequations} 
where for the sake of simplicitiy and compactness of notation, we have dropped function arguments. In Equation~\eqref{eq:optconditionslogbarrier}, $\badjoint$ is a Lagrange multiplier associated to the PDE-constraint denoted in Equation~\eqref{eq:discreteequalityconstraint}, $\bs{M}_{\! L}=\text{diag}(\bs{M}\bs{1})$ is a diagonal lumped mass matrix, $\bs{J}_{\! \bparam}$ and  $\bs{J}_{\! \bstatevar}$ are the Jacobians of the (discretized) PDE-constraint function from Equation~\eqref{eq:discreteequalityconstraint}; furthermore, $\bs{\nabla}_{\! \bstatevar}, \bs{\nabla}_{\! \bparam}$ are the gradient operators with respect to $\bstatevar$ and $\bparam$ respectively, with which we utilize the forms of $\objfun(\statevar,\param)$ and $c(\statevar,\param,\adjoint)$ to compute quantities needed for a complete description of the optimality system, described in Equation~\eqref{eq:optconditionslogbarrier}, i.e.,
\begin{align*}
(\bs{\nabla}_{\! \bstatevar}f_{h})_{i}&=\frac{\delta}{\delta \statevar_{h}}f(\statevar_{h},\param_{h})(\phi_{i}):=
\left[\frac{\partial }{\partial \varepsilon}f(\statevar_{h}+\varepsilon \phi_{i}, \param_{h})\right]_{\varepsilon=0},\quad 1\leq i\leq \dimstate,	\\
(\bs{\nabla}_{\! \bparam}f_{h})_{i}&=\frac{\delta}{\delta \param_{h}}f(\statevar_{h},\param_{h})(\psi_{i}):=
\left[\frac{\partial }{\partial \varepsilon}f(\statevar_{h}, \param_{h}+\varepsilon\psi_{j})\right]_{\varepsilon=0},\quad 1\leq i\leq \dimparam, \\
\left(\bs{J}_{\! \bstatevar}\right)_{i,j}&=
\frac{\delta}{\delta \statevar_{h}}\constfun(\statevar_{h}, \param_{h}, \phi_{i})(\phi_{j}):=
\left[
\frac{\partial}{\partial \varepsilon}\constfun(\statevar_{h}+\varepsilon\phi_{j},\param_{h},\phi_{i})
\right]_{\varepsilon=0},\quad 1\leq i\leq \dimstate,\quad 1\leq j\leq \dimstate,\\
\left(\bs{J}_{\! \bparam}\right)_{i,j}&=
\frac{\delta}{\delta \param_{h}}\constfun(\statevar_{h}, \param_{h}, \phi_{i})(\psi_{j}):=
\left[
\frac{\partial}{\partial \varepsilon}\constfun(\statevar_{h},\param_{h}+\varepsilon \psi_{j},\phi_{i})
\right]_{\varepsilon=0},\quad 1\leq i\leq \dimstate,\quad 1\leq j\leq \dimparam.\\
\end{align*}
The term $\bs{M}_{\! L}(\mu\bs{1}/(\bparam-\bparamlbound))$, is due to the gradient of the discretized log-barrier $\mu\bs{1}^{\top}\bs{M}\log(\bparam-\bparamlbound)$ and the operator $\delta/\delta \statevar_{h}$ is the variational derivative operator~\cite{gelfand1961} with respect to $\statevar_{h}$.

We formally introduce the dual variable $\blboundmultiplier=\mu\bs{1}/(\bparam-\bparamlbound)$, associated to the bound constraint $\bparam\geq\bparam_{\ell}$. With this formulation, the first-order optimality conditions require that a primal-dual solution $(\bstatevar^{\star,\mu}, \bparam^{\star,\mu}, \badjoint^{\star,\mu}, \blboundmultiplier^{\star,\mu})$ of the log-barrier subproblem~\eqref{eq:finitedimlogbarrier} satisfies the nonlinear system of equations
\begin{subequations} \label{eq:smoothKKT} 
	\begin{align}
	\bs{\nabla}_{\! \bstatevar}f_{h}+\bs{J}_{\! \bstatevar}^{\top}\badjoint^{\star,\mu}&=\bs{0},
	\label{eq:smoothKKT1}\\
	\bs{\nabla}_{\! \bparam}f_{h}+\bs{J}_{\! \bparam}^{\top}\badjoint^{\star}-\bs{M}_{\! L} \blboundmultiplier^{\star,\mu}&=\bs{0}
	, \label{eq:smoothKKT2}\\
	\bs{\constfun}&=\bs{0}
	, \label{eq:smoothKKT3} \\
	\blboundmultiplier^{\star,\mu}\odot(\bparam^{\star,\mu}-\bparamlbound)-\mu\bs{1}&=\bs{0}, \label{eq:smoothKKT4}
	\end{align}
\end{subequations} 
where for $\mu>0$, Equation~\eqref{eq:smoothKKT4} guarantees that $\blboundmultiplier^{\star,\mu}=\mu\bs{1}/(\bparam^{\star,\mu} - \bparamlbound)$. An alternative interpretation of Equation~\eqref{eq:smoothKKT} is that it is a perturation of the necessarily conditions for optimality
\begin{align*}
\bs{\nabla}_{\! \bstatevar}\mathcal{L}=\bs{0},\quad
\bs{\nabla}_{\! \bparam}\mathcal{L}=\bs{0},\quad
\bs{\nabla}_{\badjoint}\mathcal{L}=\bs{0},\quad
\blboundmultiplier^{\star,\mu}\odot(\bparam^{\star,\mu}-\bparam_{\ell})-\mu\bs{1}=\bs{0},
\end{align*}
of the modified Lagrangian
\begin{align} \label{eq:Lagrangian}
\mathcal{L}(\bstatevar,\bparam,\badjoint,\blboundmultiplier)&:=
\objfun_{h}(\bstatevar,\bparam)+\badjoint^{\top}\bs{c}(\bstatevar,\bparam)-\blboundmultiplier^{\top}\bs{M}_{\! L}(\bparam-\bparam_{\ell}),
\end{align}
associated to the discretized constrained optimization problem
\begin{subequations} 
	\begin{align}
	\min_{(\bstatevar,\bparam)\in\mathbb{R}^{\dimstate}\times\mathbb{R}^{\dimparam}}&\objfun_{h}(\bstatevar,\bparam), \\
	\text{ such that }&\bs{c}(\bstatevar,\bparam)=\bs{0},\\
	\text{ and }&\bparam\geq\bparam_{\ell}.
	\end{align}
\end{subequations}

\subsection{Stopping criteria}\label{subsec:stoppingcriteria}

Next, various metrics are presented that we use to measure progress towards
the solution of a discretization of the optimization problem~\eqref{eq:infdimproblem}, as well as a discretized log-barrier subproblem~\eqref{eq:finitedimlogbarrier}. As in~\cite{petra2023}, we utilize norms that measure proximity towards a local optimum that are well-defined with respect to an infinite-dimensional problem formulation. Such norms are used in order to avoid discretization dependent performance of the optimization algorithm. To measure stationarity of the Lagrangian (Equation~\eqref{eq:Lagrangian}) we use
\begin{align*}
e^{rr}_{\text{stat}}(\bstatevar, \bparam,\badjoint, \blboundmultiplier):=\sqrt{\|\bs{\nabla}_{\! \bstatevar}\mathcal{L}\|_{\bs{M}_{\bstatevar}^{-1}}^{2}
	+\|\bs{\nabla}_{\! \bparam}\mathcal{L}\|_{\bs{M}^{-1}}^{2}},
\end{align*}
and 
\begin{align*}
e^{rr}_{\text{feas}}(\bstatevar,\bparam):=\|\bs{c}(\bstatevar,\bparam)\|_{\bs{M}_{\bstatevar}^{-1}},
\end{align*}
to measure feasibility, where $\bs{M}_{\bstatevar}$ is the mass matrix with respect to the space $\mathcal{V}_{h}$, $(\bs{M}_{\bstatevar})_{i,j}=\int_{\Omega}\phi_{i}(\bs{x})\,\phi_{j}(\bs{x})\,\mathrm{d}\bs{x}$ and as before $\bs{M}$ is the mass matrix with respect to the space $\mathcal{M}_{h}$. To measure complementarity we make use of
\begin{align*}
e^{rr}_{\text{compl}}(\bparam, \blboundmultiplier;\mu):=\bs{1}^{\top}\bs{M}|\blboundmultiplier\odot(\bparam-\bparamlbound)-\mu\bs{1}|,
\end{align*}
where the absolute value $|\cdot|$ is taken componentwise, i.e., $|\bs{x}|_{i}=|\bs{x}_{i}|$.
The stationarity, feasibility and complementarity measures are combined into the single optimality measure
\begin{align*}
e^{rr}(\bstatevar,\bparam,\badjoint, \blboundmultiplier;\mu):=\max\left\lbrace e^{rr}_{\text{stat}}(\bstatevar, \bparam, \badjoint, \blboundmultiplier)/s_{d},
e^{rr}_{\text{feas}}(\bstatevar, \bparam), e^{rr}_{\text{compl}}(\bparam, \blboundmultiplier; \mu)/s_{c}\right\rbrace,
\end{align*}
where the scaling values $s_{c}$ and $s_{d}$ are given by
\begin{align*}
s_{c}=\max\left\lbrace
s_{\text{max}}, \|\blboundmultiplier\|_{\bs{M}}
\right\rbrace/s_{\text{max}}, \\
s_{d}=\max\left\lbrace
\frac{1}{2}\|\badjoint\|_{\bs{M}_{\bstatevar}}+\frac{1}{2}\|\blboundmultiplier\|_{\bs{M}},s_{\text{max}}
\right\rbrace/s_{\text{max}},
\end{align*}
where $s_{\text{max}}$ is a constant as in~\cite{wachter2006}. The scaling is used to avoid difficulties in reducing the optimality error when the Lagrange multipliers grow too large. We note that $e^{rr}(\bstatevar,\bparam,\badjoint,\blboundmultiplier;\mu)$ is an optimality measure for the discretized log-barrier subproblem~\eqref{eq:finitedimlogbarrier} and $e^{rr}(\bstatevar,\bparam,\badjoint,\blboundmultiplier;0)$ is an optimality measure for the discretized PDE- and bound-constrained optimization problem~\eqref{eq:infdimproblem}.

\section{An IP-Gauss-Newton method}
\label{sec:IPGaussNewton}

The nonlinear primal-dual optimality system associated to the log-barrier subproblem~\eqref{eq:smoothKKT} is solved by a damped (\textit{e.g., line-search}) Gauss-Newton method. This is achieved by the following modification of the linearization of Equation~\eqref{eq:smoothKKT}
\begin{align} \label{eq:linearizedsmoothKKT}
\begin{bmatrix}
\bs{H}_{\! \bstatevar, \bstatevar} & \bs{0} & \bs{J}_{\! \bstatevar}^{\top} & \bs{0} \\
\bs{0} & \bs{H}_{\! \bparam, \bparam} & \bs{J}_{\! \bparam}^{\top} & -\bs{M}_{\! L} \\
\bs{J}_{\! \bstatevar} & \bs{J}_{\! \bparam} & \bs{0} & \bs{0} \\
\bs{0} & \text{diag}(\blboundmultiplier) & \bs{0} & \text{diag}(\bparam-\bparamlbound)
\end{bmatrix}
\begin{bmatrix}
\bstatevarupdate \\
\bparamupdate \\
\badjointupdate \\
\blboundmultiplierupdate
\end{bmatrix}
&=
-
\begin{bmatrix}
\bs{r}_{\bstatevar} \\
\bs{r}_{\bparam} \\
\bs{r}_{\badjoint} \\
\bs{r}_{\blboundmultiplier}
\end{bmatrix},
\end{align}
where
\begin{align}
\bs{H}_{\! \bstatevar, \bstatevar}&:=\bs{\nabla}_{\! \bstatevar,\bstatevar}^{2}(
f_{h}),\quad (\bs{H}_{\! \bstatevar, \bstatevar})_{i,j}:=
\frac{\delta^{2}}{\delta \statevar_{h}^{2}}(f(\statevar_{h}, \param_{h}))(\phi_{i},\phi_{j}), \\
\bs{H}_{\! \bparam, \bparam}&=\bs{\nabla}_{\! \bparam,\bparam}^{2}(f_{h}),\quad (\bs{H}_{\! \bparam, \bparam})_{i,j}:=
\frac{\delta^{2}}{\delta \param_{h}^{2}}
(
f(\statevar_{h}, \param_{h}))(\psi_{i},\psi_{j}). \label{eq:Hrhorho}
\end{align}
The right-hand side of the linear system described by Equation~\eqref{eq:linearizedsmoothKKT} is defined in terms of the optimality residuals
\begin{align*}
\bs{r}_{\bstatevar}&=\bs{\nabla}_{\! \bstatevar}f_{h}+\bs{J}_{\! \bstatevar}^{\top}\badjoint,
\\
\bs{r}_{\bparam}&=\bs{\nabla}_{\! \bparam}f_{h}+\bs{J}_{\! \bparam}^{\top}\badjoint-\bs{M}_{\! L} \blboundmultiplier
, \\
\bs{r}_{\badjoint}&=\bs{\constfun}, \\
\bs{r}_{\blboundmultiplier}&=\blboundmultiplier\odot(\bparam-\bparamlbound)-\mu\bs{1}.
\end{align*}
The solution $(\bstatevarupdate,\bparamupdate,\badjointupdate,\blboundmultiplierupdate)$ of Equation~\eqref{eq:linearizedsmoothKKT} is the so-called Gauss-Newton search direction and is used to update a given optimizer estimate $({\bstatevar}, {\bparam}, {\badjoint}, {\blboundmultiplier})\rightarrow({\bstatevar}^+, {\bparam}^+, {\badjoint}^+, \blboundmultiplier^+)$ based on the linear update
\begin{subequations} \label{eq:newtonupdate}  
	\begin{align}
	\begin{bmatrix}
	\bstatevar^{+} & \bparam^{+} & \badjoint^{+}
	\end{bmatrix}
	&=
	\begin{bmatrix}
	\bstatevar & \bparam & \badjoint
	\end{bmatrix}
	+
	\alpha_p
	\begin{bmatrix}
	\bstatevarupdate & \bparamupdate & \badjointupdate
	\end{bmatrix},  \label{eq:newtonupdatep} \\
	\blboundmultiplier^{+}&=\blboundmultiplier+\alpha_d\,\blboundmultiplierupdate. \label{eq:newtonupdated}
	\end{align}
\end{subequations}
The primal and dual step-lengths $\alpha_p$ and $\alpha_d$ are computed using a convergence-enforcing filter line-search algorithm~\cite{wachter2006}. 
The salient idea behind the filter line-search, in the context of the log-barrier subproblem~\eqref{eq:finitedimlogbarrier}, is that a new trial point~\eqref{eq:newtonupdate} is accepted when either the log-barrier objective~\eqref{eq:logbarobjective_discrete} or the norm of the constraint~\eqref{eq:logbareqcons_discrete} is sufficiently reduced relative to the current iterate. When this condition is not satisfied, the primal step-length $\alpha_p$ is reduced. The largest possible primal and dual step-lengths $\alpha_p,\alpha_d$ are chosen using a fraction-to-boundary rule that ensures the primal and dual variables remain in the interior of $\lbrace (\bparam,\blboundmultiplier)\in\mathbb{R}^{\dimparam}\times\mathbb{R}^{\dimparam}\text{ such that }\blboundmultiplier\geq\bs{0} \text{ and }\bparam-\bparam_{\ell}\geq\bs{0} \rbrace$.

By algebraically eliminating $\blboundmultiplierupdate$, the system of Equations~\eqref{eq:linearizedsmoothKKT} can be posed as the symmetric saddle-point linear system  
\begin{align} \label{eq:IPNewtonsys}
\underbrace{
	\begin{bmatrix}
	\bs{H}_{\! \bstatevar, \bstatevar} & \bs{0} & \bs{J}_{\! \bstatevar}^{\top}   \\
	\bs{0}  & \bs{W}_{\! \bparam, \bparam}^{\mu} & \bs{J}_{\! \bparam}^{\top}  \\
	\bs{J}_{\! \bstatevar} & \bs{J}_{\! \bparam} & \bs{0}  \\
	\end{bmatrix}
}_{\bs{A}}
\underbrace{ 
	\begin{bmatrix}
	\bstatevarupdate \\
	\bparamupdate \\
	\badjointupdate	
	\end{bmatrix}
}_{\bs{x}}
&=
\underbrace{ 
	\begin{bmatrix}
	\bs{b}_{\bstatevar} \\
	\bs{b}_{\bparam}\\
	\bs{b}_{\badjoint}	
	\end{bmatrix}
}_{\bs{b}},
\end{align}
where
\begin{align}
\bs{W}_{\! \bparam, \bparam}^{\mu}&=\bs{H}_{\! \bparam, \bparam}+\bs{H}_{\! \text{log-bar}},  \label{eq:Wrhorho} \\
\bs{H}_{\! \text{log-bar}}&= \bs{M}_{\! L}\text{diag}( \blboundmultiplier/(\bparam-\bparamlbound)), \label{eq:Hlogbar}\\
\bs{b}_{\bparam}&=-(\bs{r}_{\bparam}+\bs{M}_{\! L}\bs{r}_{\blboundmultiplier}/(\bparam-\bparamlbound)),\\
\bs{b}_{\bstatevar}&=-\bs{r}_{\bstatevar},\quad\bs{b}_{\badjoint}=-\bs{r}_{\badjoint}.
\end{align}
The superscript $\mu$, contained in the symbol $\bs{W}_{\! \bparam, \bparam}^{\mu}$ indicates that said matrix depends on the log-barrier parameter $\mu$. This dependence has been made symbolically explicit since $\bs{W}_{\!\bparam,\bparam}^{\mu}$ generally becomes ill-conditioned as $\mu$ tends to zero. $\bs{H}_{\! \text{log-bar}}$ is termed the log-barrier Hessian as when perturbed complementarity, Equation~\eqref{eq:smoothKKT4}, holds, it is precisely the Hessian of the log-barrier term $-\mu\bs{1}^{\top}\bs{M}\log(\bparam-\bparam_{\ell})$ with respect to $\bparam$. Equation~\eqref{eq:IPNewtonsys}
is hereafter referred to
as the IP-Gauss-Newton linear system. Once a solution to the IP-Gauss-Newton linear system has been found, the search direction for the bound constraint Lagrange multiplier $\blboundmultiplierupdate$, is determined, by backsubstitution on the fourth block row of Equation~\ref{eq:linearizedsmoothKKT}, as
\begin{align*} 
\blboundmultiplierupdate = -\left[\blboundmultiplier+(\blboundmultiplier\odot\bparamupdate-\mu\bs{1})/(\bparam-\bparamlbound)\right].
\end{align*}

%
%

It is common to symmetrize  IP linear systems such as~\eqref{eq:linearizedsmoothKKT}, in the form of~\eqref{eq:IPNewtonsys}, for which $\bs{L}\bs{D}\bs{L}^{\top}$ factorizations can be used and are generally more stable and efficient than $\bs{L}\bs{U}$ factorizations. Symmetrization is also generally favorable for Krylov subspace methods as it permits usage of methods with lower-memory requirements than GMRES, such as CG and MINRES. The reduction of Equation~\eqref{eq:linearizedsmoothKKT} to Equation~\eqref{eq:IPNewtonsys} provides a symmetric linear system of reduced size at a negligible computational cost since the computation of $\bs{H}_{\!\text{log-bar}}$, $\bs{b}_{\bparam}$, and $\blboundmultiplierupdate$ require only elementwise vector operations. A potential disadvantage of this approach is the deterioration of the condition number of the linear system matrix due to the diagonal log-barrier Hessian 
$\bs{H}_{\! \text{log-bar}}$. The condition number of the log-barrier Hessian, in general, grows unboundedly as $\mathcal{O}(\mu^{-2})$ as the IPM progresses towards the optimal solution~\cite{nocedalwright2006}. Other symmetrization strategies for Equation~\eqref{eq:linearizedsmoothKKT} are possible. Recent work by Ghannad et al.~\cite{GhannadOrbanSaunders_21_LinsysIpm} (also see related previous theoretical work by Greig et al.~\cite{GreigMouldingOrban_14_BoundIpm} that bound the eigenvalues and condition number of various formulations of IP linear systems) indicates that symmetrizations of Equation~\eqref{eq:linearizedsmoothKKT} that keep the $4\times 4$ structure are more well conditioned than the reduced symmetric linear system~\eqref{eq:IPNewtonsys} for convex problems with proper regularizations. This likely applies to our nonconvex setup; in fact, the state-of-the-art solver Ipopt works with one such $4\times 4$ symmetric linear system~\cite{wachter2006} . However, in the case of the problems considered in this work, the ill-conditioning of the symmetric $3\times 3$ linear system~\eqref{eq:IPNewtonsys} can be effectively factored out by preconditioning strategies such as those elaborated in Section~\ref{sec:linearalgebra}, making the size reduction worthy.

A notable difference of our approach from the filter line-search algorithm presented in~\cite{wachter2005} is the use of a Gauss-Newton search direction instead of a Newton search direction with inertia regularization. Inertia regularization is a mechanism needed to ensure that the solution of the Newton linear system is a descent direction of the log-barrier objective when the equality constraint~\eqref{eq:smoothKKT3} is nearly satisfied, and is critical to ensure convergence of the line-search based IPM~\cite{wachter2005}. A sufficient condition for a Newton search direction to be of descent is that the inertia of the Newton linear system matrix is the triplet $(\dimstate+\dimparam, \dimstate, 0)$. However, it is not feasible to compute the inertia of the Newton linear system matrix since we target large-scale problems and utilize Krylov subspace methods for approximate linear solves. Another inertia regularization technique is inertia-free regularization~\cite{chiang2016}, that repeatedly perturbs diagonal blocks of the Newton linear system matrix until carefully designed ``curvature'' tests are satisfied.
Both inertia-based and inertia-free techniques require additional Newton linear system solves as well as matrix-vector products (in addition to the modification of the diagonal subblocks of the system matrix).

A Gauss-Newton approach does not require computationally expensive linear algebraic operations needed to detect and signal inertia-correction. This is because it generates search directions that are guaranteed to be of descent for points that nearly satisfy the PDE-constraint. However, the Gauss-Newton method does not asymptotically converge at a quadratic rate like Newton's method.
The Gauss-Newton search direction being of descent when the PDE-constraint is nearly satisfied is a consequence of the upper $2\times 2$ block 
\begin{align*}
\begin{bmatrix}
\bs{H}_{\! \bstatevar, \bstatevar} & \bs{0} \\
\bs{0} & \bs{W}_{\! \bparam, \bparam}^{\mu}
\end{bmatrix},
\end{align*}
of~\eqref{eq:IPNewtonsys} being positive semi-definite. We next project said $2\times 2$ block matrix onto the nullspace of the Jacobian of the PDE-constraint with respect to $(\bstatevar,\bparam)$. The columns of 
\begin{align*}
\bs{N}=
\begin{bmatrix}
-\bs{J}_{\! \bstatevar}^{-1}\bs{J}_{\! \bparam} \\
\bs{I} 
\end{bmatrix},
\end{align*}
form a basis for the nullspace of PDE-constraint Jacobian $\begin{bmatrix} \bs{J}_{\! \bstatevar} & \bs{J}_{\! \bparam}\end{bmatrix}$ and that the Hessian of $\objfun_{h}$ with respect to $(\bstatevar,\bparam)$ projected onto said null-space is
\begin{align}
\bs{\hat{H}}&=
\bs{N}^{\top}
\begin{bmatrix}
\bs{H}_{\! \bstatevar, \bstatevar} & \bs{0} \\%
\bs{0} & \bs{W}_{\! \bparam, \bparam}^{\mu}
\end{bmatrix}
\bs{N} \nonumber \\
&=\bs{\hat{H}}_{\bs{d}}+\bs{W}_{\! \bparam, \bparam}^{\mu} \label{eq:reduced_hessian}\\
\bs{\hat{H}}_{\bs{d}}:&=
(\bs{J}_{\! \bstatevar}^{-1}\bs{J}_{\! \bparam})^{\top}\bs{H}_{\! \bstatevar, \bstatevar}
(\bs{J}_{\! \bstatevar}^{-1}\bs{J}_{\! \bparam}).\label{eq:reduced_misfithessian}	
\end{align}
$\bs{\hat{H}}$ is known as the reduced-space Gauss-Newton Hessian in the context of nonlinear least squares~\cite[Chapter 10.3]{nocedalwright2006} and its relation to $\bs{A}$
is precisely what motivates us to refer to the matrix $\bs{A}$, as defined in Equation~\eqref{eq:IPNewtonsys}, as a IP-Gauss-Newton linear system matrix. It is to be noted that, while we do not use a reduced-space approach for optimization, it is the case that typical reduced-space PDE-constrained optimization problems~\cite{ghattas2021,hartland2022thesis}, are nonlinear least squares problems.

\section{Preconditioners for scalable IP-Gauss-Newton iterative linear solvers}\label{sec:linearalgebra}
A common approach to solve indefinite saddle-point interior-point linear systems is a direct solver such as in~\cite{duff2004}. However, the computational cost increases at an undesirably high rate with respect to problem size and that limits the sizes of problems that can be solved with such approaches. PDE- and bound-constrained optimization problems is one such class of problems, wherein the dimension of the discretized optimization variable is made arbitrarily large under mesh refinement.
Here, we take advantage of problem specific structure and propose two intimately related preconditioned Krylov subspace linear system solution strategies, whose scaling relies only
on underlying multigrid linear solves of submatrix blocks that are known to be amenable to multigrid. The described framework achieves algorithmic scalability, in the sense that the number of outer optimization iterations and inner Krylov subspace iterations is largely independent of both the problem discretization and the log-barrier parameter $\mu$. Due to the aforementioned algorithmic scalability and an implementation that leverages scalable libraries, we obtain an efficient and scalable means to compute solutions of elliptic PDE- and bound-constrained optimization problems with results reported in Section~\ref{subsec:example2}.    

\subsection{A block Gauss-Seidel preconditioner for full-space IP-Gauss-Newton GMRES solves} \label{subsec:GMRESpreconditioner}

We next introduce, for Equation~\eqref{eq:IPNewtonsys}, the block Gauss-Seidel preconditioner
\begin{align} \label{eq:BlockGSPreconditioner}
\bs{\tilde{A}}&=
\begin{bmatrix} 
\bs{H}_{\! \bstatevar, \bstatevar}  & \bs{0} & \bs{J}_{\! \bstatevar}^{\top}\\
\bs{0} & \bs{W}_{\! \bparam, \bparam}^{\mu} & \bs{J}_{\! \bparam}^{\top} \\
\bs{J}_{\! \bstatevar} & \bs{0} & \bs{0} \\
\end{bmatrix}.
\end{align} 
The block Gauss-Seidel preconditioner belongs to the class of central null preconditioners~\cite{pestana2016null}, where
the Schur-complement $\bs{\hat{H}}$ is approximated by $\bs{W}_{\! \bparam, \bparam}^{\mu}$. If the Schur complement were not approximated then, in exact arithmetic, the preconditioned GMRES solve would converge in two iterations since the minimum polynomial for the preconditioned matrix is degree two~\cite{benzi2005numerical}. A number of papers advocate this and related strategies for building efficient (block triangular) preconditioners for saddle-point linear systems e.g.,~\cite{bank1989, zulehner2000, benzi2005numerical, drzisga2018, pestana2016null, pearson2018refined} and references therein. We refer to $\bs{\tilde{A}}$ as a block Gauss-Seidel preconditioner since under a suitable symmetric permutation of its rows and columns, see Equation~\eqref{eq:BlockGSupperTriangular}, it is block upper triangular. As will be shown, and demonstrated numerically in Section~\ref{sec:numericalresults}, the spectral properties of the proposed block Gauss-Seidel preconditioned IP-Gauss-Newton matrix $\bs{\tilde{A}}^{-1}\bs{A}$ has desirable properties for many large-scale PDE- and bound-constrained optimization problems. Such desirable spectral properties include asymptotic independence with respect to mesh refinement and log-barrier parameter $\mu$. 
This preconditioner is particularly effective as it exploits known structure of PDE-constrained optimization problems, namely, Property~\ref{propertyP}.
%
%
\begin{customproperty}{(P)}
	\label{propertyP}
	The eigenvalues of the regularization preconditioned reduced-space data-misfit Gauss-Newton Hessian 
	$\bs{H}_{\! \bparam, \bparam}^{-1}\bs{\hat{H}_{d}}$
	decay rapidly to zero and in a discretization independent manner~\cite{flath2011, isaac2015, ghattas2021}.
\end{customproperty}

%
%
This property stems from the limited amount of useful information~\cite[Chapter 4]{Alger19} that the data provides to the solution of the optimization problem~\eqref{eq:infdimproblem}. By useful information, we refer to the subset of data whose influence on the solution is not made negligible by the regularization component of the objective functional. This structure is reflected in the spectral properties of the regularization preconditioned reduced-space data-misfit Gauss-Newton Hessian.

The following proposition states a relation between the spectra of $\bs{\tilde{A}}^{-1}\bs{A}$ and $\bs{H}_{\! \bparam, \bparam}^{-1}\bs{\hat{H}}_{\bs{d}}$, and demonstrates that under Property~\ref{propertyP} the eigenvalues of $\bs{\tilde{A}}^{-1}\bs{A}$ rapidly decays to $1$ at a rate that is asymptotically independent of discretization and ill-conditioning coming from the log-barrier parameter $\mu$. 


\begin{proposition}\label{prop:eigs}
	Let $\bs{A}$ and $\bs{\tilde{A}}$ be specified by Equation~\eqref{eq:IPNewtonsys} and Equation~\eqref{eq:BlockGSPreconditioner} respectively. If the subblock $\bs{H}_{\bstatevar,\bstatevar}$ is positive semidefinite, the subblock $\bs{W}_{\! \bparam, \bparam}^{\mu}$ given by Equation~\eqref{eq:Wrhorho} is positive definite, and  $(\bs{W}_{\! \bparam, \bparam}^{\mu}-\bs{H}_{\! \bparam, \bparam})$ is positive semidefinite, Then, the eigenvalues of $\bs{\tilde{A}}^{-1}\bs{A}$ satisfy
	\begin{align}\label{eq:eigsarebounded}
	1 \leq \lambda_{j}(\bs{\tilde{A}}^{-1}\bs{A}) \leq
	\begin{cases}
	1+\lambda_{j}(\bs{H}_{\! \bparam, \bparam}^{-1}\bs{\hat{H}_{d}}),\quad &1\leq j\leq \dimparam,\\
	1,\quad &\dimparam+1\leq j \leq \dimparam+2\dimstate.
	\end{cases}
	\end{align}
	%
	%
	%
	\begin{proof}
		It is shown in~\cite[Theorem 2.5]{pestana2016null} that an eigenvalue of $\bs{\tilde{A}}^{-1}\bs{A}$ is either equal to $1$ or to an eigenvalue of $(\bs{W}_{\! \bparam, \bparam}^{\mu})^{-1}\bs{\hat{H}}=\bs{I}+(\bs{W}_{\! \bparam, \bparam}^{\mu})^{-1}\bs{\hat{H}}_{\bs{d}}$. Then, with respect to a descending order of the eigenvalues we have
		\begin{align} \label{eq:eigRelation}
		\lambda_{j}(\bs{\tilde{A}}^{-1}\bs{A}) =
		\begin{cases}
		1+\lambda_{j}((\bs{W}_{\! \bparam, \bparam}^{\mu})^{-1}\bs{\hat{H}_{d}}),\quad &1\leq j\leq \dimparam,\\
		1,\quad &\dimparam+1\leq j \leq \dimparam+2\dimstate.
		\end{cases}
		\end{align}
		Since $\bs{W}_{\! \bparam, \bparam}^{\mu}$ and $\bs{H}_{\! \bparam, \bparam}$ are symmetric positive definite with $\bs{W}_{\! \bparam, \bparam}^{\mu}\geq \bs{H}_{\! \bparam, \bparam}$ and $\bs{\hat{H}}_{\bs{d}}$ is positive semi-definite
		then a straightforward application of Proposition~\ref{prop:eigenvalueOrdering}, as detailed in Appendix~\ref{sec:eigOrdering}, allows one to conclude that $
		0\leq \lambda_{j}((\bs{W}_{\! \bparam, \bparam}^{\mu})^{-1}\bs{\hat{H}}_{\bs{d}})
		\leq 
		\lambda_{j}(\bs{H}_{\! \bparam, \bparam}^{-1}\bs{\hat{H}}_{\bs{d}})$.	
	\end{proof}
\end{proposition}

If, as is often the case, $\bs{H}_{\!\bstatevar,\bstatevar}$ is rank-deficient, then $\bs{\tilde{A}}^{-1}\bs{A}$ is generally not diagonalizable, see Note~\ref{note:notDiagnolizable} in the Appendix. This rank deficiency is due to the functional forms that $\objfun_{\text{misfit}}(\statevar)$ has in various PDE-constrained optimization problems. For instance, when $\objfun_{\text{misfit}}(\statevar)$ only depends on the evaluation of the state at a small discretization independent number $k$, of observation points $\objfun_{\text{misfit}}(\statevar)=\frac{1}{2}\sum_{i=1}^{k}(u(\spatialvar_{i})-\bs{d}_{i})^{2}$, then $\bs{H}_{\!\bstatevar,\bstatevar}=\bs{\mathcal{B}}^{\top}\bs{\mathcal{B}}$, where $\bs{\mathcal{B}}\in\mathbb{R}^{k\times \dimstate}$ is a so-called observation operator~\cite{ghattas2021}, such that $(\bs{\mathcal{B}}\bstatevar)_{i}\approx\statevar(\spatialvar_{i})$. For this case, $\text{Rank}(\bs{H}_{\!\bstatevar,\bstatevar})\leq k$, which is generally less than $\dimstate$. One remedy, to ensure diagonalizability, is to make a small modification to the IP-Gauss-Newton matrix as in Appendix~\ref{sec:perturbedIPGaussNewton}, for which the perturbed block Gauss-Seidel preconditioned IP-Gauss-Newton matrix is diagonalizable. The convergence behavior of the associated preconditioned GMRES solve is amenable to analysis by~\cite[Proposition 4]{saad1986}. The convergence rate of the perturbed approach is shown in Proposition~\ref{prop:GMRESresidualReduction}, in the Appendix. As detailed in Remark~\ref{rmk:perturbation} and Remark~\ref{rmk:perturbedresiduals} in the Appendix, the GMRES iteration complexity can be independent of the discretization and log-barrier parameter with the perturbed approach. We have been unable to prove but have observed in numerous numerical experiments, see Table~\ref{tbl:algorithmicscaling}, that the unperturbed Gauss-Seidel preconditioned GMRES solves converges with a number of iterations that is mesh and log-barrier independent. This suggests that the convergence of the GMRES solves for PDE- and bound-constrained problems is largely due to the spectrum of the block Gauss-Seidel preconditioned IP-Gauss-Newton matrix. 


Among alternatives to the block Gauss-Seidel preconditioner~\eqref{eq:BlockGSPreconditioner} we mention
\begin{align*}
\bs{\tilde{A}}_{\text{cen}}&=
\begin{bmatrix}
\bs{H}_{\! \bstatevar, \bstatevar} & \bs{0} & \bs{J}_{\! \bstatevar}^{\top} \\
\bs{0} & \bs{W}_{\! \bparam, \bparam}^{\mu} & \bs{0} \\
\bs{J}_{\! \bstatevar} & \bs{0} & \bs{0} \\ 
\end{bmatrix},\text{ and }
\bs{\tilde{A}}_{\text{con}}=
\begin{bmatrix}
\bs{H}_{\! \bstatevar, \bstatevar} & \bs{0} & \bs{J}_{\! \bstatevar}^{\top}  \\
\bs{0} & \bs{W}_{\! \bparam, \bparam}^{\mu} & \bs{J}_{\! \bparam}^{\top} \\
\bs{J}_{\! \bstatevar} & \bs{0}  & \bs{0} 
\end{bmatrix}
\begin{bmatrix}
\bs{I} & \bs{J}_{\! \bstatevar}^{-1}\bs{J}_{\! \bparam} & \bs{0} \\
\bs{0} & \bs{I} & \bs{0} \\
\bs{0} & -\bs{J}_{\! \bstatevar}^{-\top}\bs{H}_{\! \bstatevar, \bstatevar}\bs{J}_{\! \bstatevar}^{-1}\bs{J}_{\! \bparam} & \bs{I}
\end{bmatrix},
\end{align*}
known as central and constraint preconditioners~\cite{pestana2016null}, respectively. Applying $\bs{\tilde{A}}_{\text{cen}}^{-1}$, and $\bs{\tilde{A}}_{\text{con}}^{-1}$ to vectors 
both require a similar sequence of block solves as what is required to apply $\bs{\tilde{A}}^{-1}$ to a vector. However, the eigenvalue distribution of $\bs{\tilde{A}}_{\text{cen}}^{-1}\bs{A}$ is not as favorable as that of the block Gauss-Seidel preconditioned IP-Newton system matrix. As to be expected, we found (see Figure~\ref{fig:GMRESMINRESitsVsMu}) $\bs{\tilde{A}}_{\text{cen}}$ to be a less effective preconditioner than the proposed preconditioner $\bs{\tilde{A}}$ for IP-Gauss-Newton GMRES solves on PDE- and bound-constrained optimization problems. As
discussed in~\cite[Theorem 2.8]{pestana2016null}, the constraint preconditioner has an eigenvalue distribution equal to $(\bs{W}_{\! \bparam, \bparam}^{\mu})^{-1}\bs{\hat{H}}$. However, the application of $\bs{\tilde{A}}_{\text{con}}$ requires an additional $\bs{J}_{\! \bstatevar}$ solve and $\bs{J}_{\! \bstatevar}^{\top}$ per preconditioner apply thus making it roughly twice as expensive than $\bs{\tilde{A}}$ per apply. While both the central and constraint preconditioners are symmetric, they are
guaranteed to be indefinite thus
making them unsuitable for symmetric Krylov subspace solvers such as MINRES~\cite{van2003}, that requires symmetric positive definiteness of a preconditioner. Yet another symmetric block Gauss-Seidel preconditioner is
\begin{align*} 
\bs{\tilde{A}}_{2}&=
\begin{bmatrix}
\bs{0} & \bs{0} & \bs{J}_{\! \bstatevar}^{\top} \\
\bs{0} & \bs{W}_{\! \bparam, \bparam}^{\mu} & \bs{J}_{\! \bparam}^{\top} \\
\bs{J}_{\! \bstatevar} & \bs{J}_{\! \bparam} & \bs{0} \\ 
\end{bmatrix},
\end{align*}
that is presented in~\cite[Equation 3.6]{birosghattas2005}. We have chosen not to explore $\bs{\tilde{A}}_{2}$ in this work, as it is expected to not offer any computational benefits with respect to the block Gauss-Seidel preconditioner. Computational benefits are not expected since this preconditioner clusters eigenvalues of the IP-Gauss-Newton system identical to that of the block Gauss-Seidel preconditioner. Also, solves with the same three submatrices are required to apply this preconditioner as in the application of the block Gauss-Seidel preconditioner.

\subsection{A $\bs{W}_{\! \bparam, \bparam}^{\mu}$ preconditioner for reduced-space CG solves} \label{subsec:CGpreconditioner}
One means to use an iterative solver that exploits symmetry is to first reformulate $\bs{A}\bs{x}=\bs{b}$, by taking the Schur complement with respect to $\bparam$, that is eliminating $\bs{x}_{\bstatevar}$ and $\bs{x}_{\badjoint}$ from Equation~\eqref{eq:IPNewtonsys},
\begin{align*}
\bs{x}_{\bstatevar}&=\bs{J}_{\! \bstatevar}^{-1}(\bs{b}_{\badjoint}-\bs{J}_{\! \bparam}\bs{x}_{\bparam}),\\
\bs{x}_{\badjoint}&=\bs{J}_{\! \bstatevar}^{-\top}(\bs{b}_{\bstatevar}-\bs{H}_{\! \bstatevar, \bstatevar}\bstatevar).
\end{align*}
The resultant linear system
\begin{align} \label{eq:reduced_hessian_system}
\bs{\hat{H}}\bs{x}_{\bparam}=\bs{\hat{b}},
\end{align}
is in the unknown $\bs{x}_{\bparam}$ and the system matrix is the Schur-complement $\bs{\hat{H}}$ (see Equation~\eqref{eq:reduced_hessian}),
also known as the reduced-space Gauss-Newton Hessian. Furthermore, $\bs{\hat{b}}=\bs{b}_{\bparam}-\bs{J}_{\! \bparam}^{\top}\bs{J}_{\! \bstatevar}^{-\top}(\bs{b}_{\bstatevar}-\bs{H}_{\! \bstatevar, \bstatevar}\bs{J}_{\! \bstatevar}^{-1}\bs{b}_{\badjoint})$, is the associated reduced right hand side.
The reduced-space Gauss-Newton Hessian $\bs{\hat{H}}$ is symmetric and guaranteed to be positive definite for problems, such as that studied here, where $\bs{W}_{\! \bparam, \bparam}^{\mu}$ is positive definite and $\bs{H}_{\! \bstatevar, \bstatevar}$ is positive semi-definite. Thus, one can solve the reduced-space Gauss-Newton system $\bs{\hat{H}}\bs{x}_{\bparam}=\bs{\hat{b}}$ with a $\bs{W}_{\! \bparam, \bparam}^{\mu}$ preconditioned CG solver.

\subsection{Computational costs of the proposed preconditioners} \label{subsec:preconditionerCost}

Having established that the spectra of the proposed preconditioned linear systems, discussed in Section~\ref{subsec:GMRESpreconditioner} and Section~\ref{subsec:CGpreconditioner} respectively, are independent of discretization and the log-barrier parameter, we now discuss a scalable means to apply the preconditioners and the reduced-space data-misfit Gauss-Newton Hessian $\bs{\hat{H}}_{\bs{d}}$~\eqref{eq:reduced_hessian}. A scalable means to apply these linear operators enables the scalable solution of the elliptic PDE- and bound-constrained optimization example problem detailed in Section~\ref{sec:problemdescription} and as numerically demonstrated in Figure~\ref{fig:scaling}.   

As expected of a block triangular matrix solve and as described in Algorithm~\ref{alg:blockGSapply}, applying $\bs{\tilde{A}}^{-1}$ to a vector requires a sequence of solves with subblocks of $\bs{\tilde{A}}$. Applying the reduced-space Gauss-Newton data-misfit Hessian $\bs{\hat{H}}_{\bs{d}}$, to a vector
requires two block solves, as detailed in Algorithm~\ref{alg:HessianAction}. 
Thus the solution of three linear systems is required per $\bs{W}_{\! \bparam, \bparam}^{\mu}$ preconditioned Krylov iteration for the reduced-space Gauss-Newton Hessian, with system matrices identical to that required for the application of the block Gauss-Seidel preconditioner $\bs{\tilde{A}}$. Thus the cost to apply $\bs{\tilde{A}}^{-1}$ to a vector as well as applying $(\bs{W}_{\! \bparam, \bparam}^{\mu})^{-1}\bs{\hat{H}}$ a vector critically depends on the cost to apply $\bs{J}_{\! \bstatevar}^{-1}$, $\bs{J}_{\! \bstatevar}^{-\top}$ and $(\bs{W}_{\! \bparam, \bparam}^{\mu})^{-1}$.
The Jacobian of the PDE-constraint with respect to $\bstatevar$, $\bs{J}_{\! \bstatevar}$, is amenable to an algebraic multigrid treatment, e.g., when the PDE-constraint~\eqref{eq:discreteequalityconstraint} describes a discretized elliptic PDE. $\bs{W}_{\! \bparam, \bparam}^{\mu}$ is the sum of $\bs{H}_{\! \bparam, \bparam}$, the Hessian of the objective $\objfun$ with respect to the parameter and the log-barrier Hessian (see Equation~\eqref{eq:Hlogbar}). 
$\bs{H}_{\! \bparam, \bparam}$ is the Hessian of the regularization term of the objective with respect to $\bparam$. For the example problem described in~\ref{sec:problemdescription}, the regularization $\objfun_{\text{reg}}$ is a linear combination of the squared $L^{2}(\Omega)$ norm of $\param$ and the squared $L^{2}(\Omega)$ norm of the gradient of $\param$, for this reason the Hessian of the regularization is a linear combination of mass and stiffness matrices, hence an invertible discretized second order elliptic PDE operator. The elliptic structure of $\bs{H}_{\! \bparam, \bparam}$~\eqref{eq:Hrhorho} is one reason why AMG preconditioned CG is a performant and scalable means to solve linear systems with a $\bs{W}_{\! \bparam, \bparam}^{\mu}$ system matrix.  Remarkably, while the ill-conditioning of the log-barrier Hessian $\bs{H}_{\! \text{log-bar}}$ (Equation~\eqref{eq:Hlogbar}) negatively impacts the performance of many
Krylov subspace based strategies for the IP-(Gauss-)Newton system~\eqref{eq:IPNewtonsys}, here the convergence of an algebraic multigrid preconditioned conjugate gradient solve (CG-AMG) of $\bs{W}_{\! \bparam, \bparam}^{\mu}$ is actually accelerated by the (diagonal) positive-definite Hessian of the log-barrier term due to improved diagonal dominance that is a key component of many AMG smoothers. Furthermore, as shown in Section~\ref{subsec:GMRESpreconditioner} the condition number of $\bs{\tilde{A}}^{-1}\bs{A}$ is independent of the log-barrier Hessian. 
\begin{algorithm} 
	\caption{
		Application of the block Gauss-Seidel preconditioner $\bs{x}=\bs{\tilde{A}}^{-1}\bs{b}$, where $\bs{x}=\begin{bmatrix}\bs{x}_{\bstatevar} & \bs{x}_{\bparam} & \bs{x}_{\badjoint}\end{bmatrix}^{\top}$ and $\bs{b}=\begin{bmatrix}
		\bs{b}_{\bstatevar} & \bs{b}_{\bparam} & \bs{b}_{\badjoint}\end{bmatrix}^{\top}$.}\label{alg:blockGSapply}
	\begin{algorithmic} 
		\State{Compute $\bs{x}_{\bstatevar}=\bs{J}_{\! \bstatevar}^{-1}\bs{b}_{\badjoint}$\hfill $\lbrace$AMG preconditioned CG solve$\rbrace$}
		\State{Compute  $\bs{x}_{\badjoint}=\bs{J}_{\! \bstatevar}^{-\top}(\bs{b}_{\bstatevar}-\bs{H}_{\! \bstatevar, \bstatevar}\bs{x}_{\bstatevar})$\hfill $\lbrace$AMG preconditioned CG solve$\rbrace$}
		\State{Compute $\bs{x}_{\bparam}=(\bs{W}_{\! \bparam, \bparam}^{\mu})^{-1}(\bs{b}_{\bparam}-\bs{J}_{\! \bparam}^{\top}\bs{x}_{\badjoint})$\hfill $\lbrace$AMG preconditioned CG solve$\rbrace$\phantom{.}}
	\end{algorithmic}
\end{algorithm}

\begin{algorithm}
	\caption{Reduced-space Gauss-Newton data-misfit Hessian-vector product $\bs{y} =\bs{\hat{H}}_{\bs{d}}\bs{x}$.} \label{alg:HessianAction}
	\begin{algorithmic}
		\State{Compute $\bs{x}_{\bstatevar}=-\bs{J}_{\! \bstatevar}^{-1}\left(\bs{J}_{\! \bparam}\bs{x}\right)$\hfill $\lbrace$AMG preconditioned CG solve$\rbrace$}
		\State{Compute $\bs{x}_{\badjoint}=-\bs{J}_{\! \bstatevar}^{-\top}\left(\bs{H}_{\! \bstatevar, \bstatevar}\bs{x}_{\bstatevar}\right)$\hfill $\lbrace$AMG preconditioned CG solve$\rbrace$}
		\State{Compute $\bs{y}\phantom{e}=\bs{J}_{\! \bparam}^{\top}\bs{x}_{\badjoint}$\hfill $\lbrace$Matrix-vector product$\rbrace$\phantom{.}}
	\end{algorithmic}	
\end{algorithm}

\section{Numerical results}\label{sec:numericalresults}

In this section, we present numerical results of the proposed framework. The results were generated via a C++ implementation of the computational framework and is a performant and distributed memory parallel means of numerically estimating optimizers of discretized PDE- and bound-constrained optimization problems. The implementation makes extensive use of the modular finite elements library MFEM~\cite{andrej2024} for finite-element discretization and also uses scalable AMG technology from hypre~\cite{hypre}. In Section~\ref{sec:problemdescription} we present a nonlinear elliptic PDE- and bound-constrained optimization example problem. In Section~\ref{subsec:example2} we show the solution of the example problem along with details of the algorithmic scaling of the proposed framework with respect to the number of IP-Gauss-Newton steps and the number of preconditioned Krylov subspace iterations for the IP-Gauss-Newton linear system solves. Numerical evidence is provided that supports Proposition~\ref{prop:GMRESresidualReduction}, namely that the performance of the proposed preconditioners is asymptotically independent of the log-barrier parameter $\mu$. We further show, having utilized mature and performant finite element~\cite{andrej2024} and algebraic multigrid~\cite{hypre} software libraries, favorable parallel scalability of the solution computation of the nonlinear PDE- and bound-constrained example problem. In Section~\ref{subsec:Morozov} we detail how the regularization parameters, $\gamma_{1},\gamma_{2}$, that enter the example problem in Equation~\eqref{eq:regularization}, are chosen. Table~\ref{tbl:noiseTable} and Table~\ref{tbl:regRobustnessTable} highlights the regularization parameter dependence of the number of GMRES iterations needed to solve Gauss-Seidel preconditioned IP-Gauss-Newton linear systems.

\subsection{Problem setup}\label{sec:problemdescription}
Here, we detail the PDE- and bound-constrained optimization problem~\eqref{eq:infdimproblem}, to determine functions $\statevar(\spatialvar)$, $\param(\spatialvar)$ defined over the closure of the unit square $\Omega=(0,1)\times(0,1)$, that is used to test the computational performance of the framework detailed in this work and from which we obtain the results that are presented in Section~\ref{sec:numericalresults}. The objective functional, whose minimizer is saught, is
a linear combination of so-called data-misfit $\objfun_{\text{misfit}}$ and regularization $\objfun_{\text{reg}}$ terms
\begin{align*}
\objfun(\statevar,\param)=\objfun_{\text{misfit}}(\statevar)+\objfun_{\text{reg}}(\param).
\end{align*}
Where,
\begin{align*}
\objfun_{\text{misfit}}(\statevar)=\frac{1}{2}\int_{\Omega_{\text{left}}}(\statevar(\spatialvar)-\statevar_{d,\zeta}(\spatialvar))^{2}\mathrm{d}\spatialvar,
\end{align*}
provides a measure of discrepancy
between the state $\statevar$ and noisy data $\statevar_{d,\zeta}$ over the left hand side $\Omega_{\text{left}}=(0,0.5)\times(0,1)$ of the domain $\Omega$. The symbols $d$ and $\zeta$ contained in the subscript of $\statevar_{d,\zeta}$ respectfully refer to data and noise. The regularization term 
\begin{align} \label{eq:regularization}
\objfun_{\text{reg}}(\param)=
\frac{\gamma_{1}}{2}
\int_{\Omega}\param^{2}(\spatialvar)\,\mathrm{d}\spatialvar+
\frac{\gamma_{2}}{2}\int_{\Omega}\bs{\nabla}_{\! \spatialvar}\param(\spatialvar)\cdot\bs{\nabla}_{\! \spatialvar}\param(\spatialvar)\,\mathrm{d}\spatialvar,
\end{align}
penalizes the squared $L^{2}(\Omega)$ norms of $\param(\spatialvar)$ and $\bs{\nabla}_{\! \spatialvar}\param(\spatialvar)$, and reduces the sensitivity of the solution $(\statevar^{\star},\param^{\star})$ of the optimization problem~\eqref{eq:infdimproblem} to random noise $\zeta$, contained in the data $u_{d,\zeta}$. The partial differential equality constraint, $c$, is the weak form: find $u\in H^{1}(\Omega)$ such that
\begin{align} \label{eq:examplePDEconstraint}
\int_{\Omega}\left(\rho\bs{\nabla}_{\! \spatialvar}\statevar\cdot\bs{\nabla}_{\! \spatialvar}\adjoint 
+\adjoint(\statevar+\statevar^{3}/3-g)\right)\mathrm{d}\spatialvar=0,\,\,\,\forall \adjoint\in H^{1}(\Omega).
\end{align}
When this elliptic PDE-constraint is expressed in strong form it reads
\begin{align*}
-\bs{\nabla}_{\! \spatialvar}\cdot(\param \bs{\nabla}_{\! \spatialvar}\statevar)+\statevar+\statevar^{3}/3 &=g,\quad\text{ in }\Omega, \\
\param \bs{\nabla}_{\! \spatialvar}\statevar\cdot\bs{n}&=0,\quad \text{ on }\partial\Omega,
\end{align*}
where $\partial\Omega$ is the boundary of the spatial domain $\Omega$, $\bs{n}$ is the outward normal to $\partial\Omega$ and $g=g(\spatialvar)$ is a forcing term. The forcing term $g=-\bs{\nabla}_{\! \spatialvar}\cdot\left(\param_{\text{true}}\,\bs{\nabla}_{\! \spatialvar}\statevar_{d}\right)+\statevar_{d}+\statevar_{d}^{3}/3$, is chosen so that the noise free data $\statevar_{d}(\spatialvar)=\cos(\pi\,\spatialvar_{1})\,\cos(\pi\,\spatialvar_{2})$ solves the PDE for the true parameter field $\param_{\text{true}}(\spatialvar)=1.0+\spatialvar_{2}\,\exp(-\spatialvar_{1}^{2})$. 

The noisy data $\statevar_{d,\zeta}$ is synthetically generated by adding random noise $\zeta$ to $\statevar_{d}$. 
The noise is obtained by taking $\zeta$ to be a random sample from a zero mean Gaussian distribution with a bi-Laplacian covariance matrix, $(-\gamma_{\zeta}\Delta+\delta_{\zeta})^{-2}$, i.e., a covariance matrix that is the square of a discretized inverse elliptic PDE operator see e.g.,~\cite{villa2021}. The motivation behind choosing this covariance is that samples from this distribution are well defined with respect to mesh refinement. A noise sample, such as that shown in Figure~\ref{fig:Morozov} (left), is uniformly scaled so that $\|\zeta\|_{L^{2}(\Omega)}=\sigma_{\zeta}\|\statevar_{d}\|_{L^{2}(\Omega)}$ for a given noise level, $\sigma_{\zeta}$. The correlation length of the noise is chosen as $0.25$, and the noise level is $5\%$ unless otherwise indicated.  While $\param_{\text{true}}$ does satisfy the lower-bound constraint $\rho\geq \param_{\ell}=1.0$, the optimizer $\param^{\star}$ can be active on a subset of $\Omega$, due to $\objfun_{\text{reg}}$ promoting small values of $\param$ when $\gamma_{1}> 0$. The constants $\gamma_{1},\gamma_{2}$ that define the regularization functional $\objfun_{\text{reg}}(\param)$ are set equal to one another and their shared value is then chosen according to the Morozov discrepancy principle, see Section~\ref{subsec:Morozov} for more details.


\subsection{Parallel and algorithmic scaling of the IPM framework}\label{subsec:example2}

\begin{figure}[h!] 
	\begin{center}
		\begin{tabular}[c]{ccc}
			\includegraphics[scale=.15]{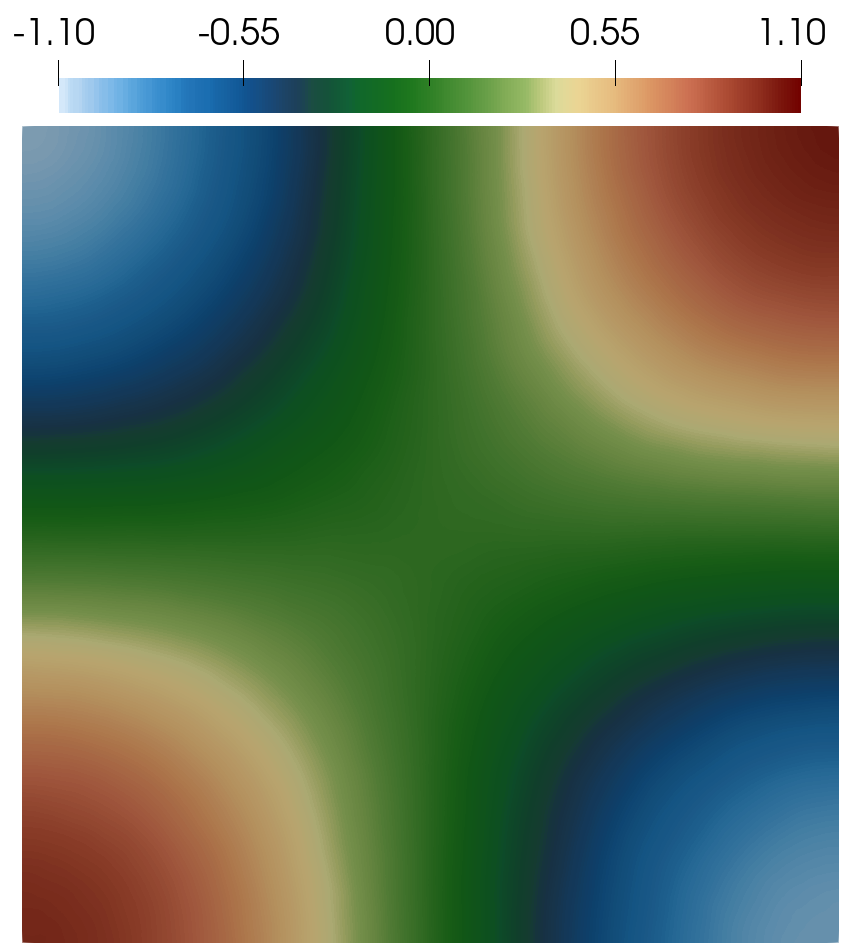} & 
			\includegraphics[scale=.15]{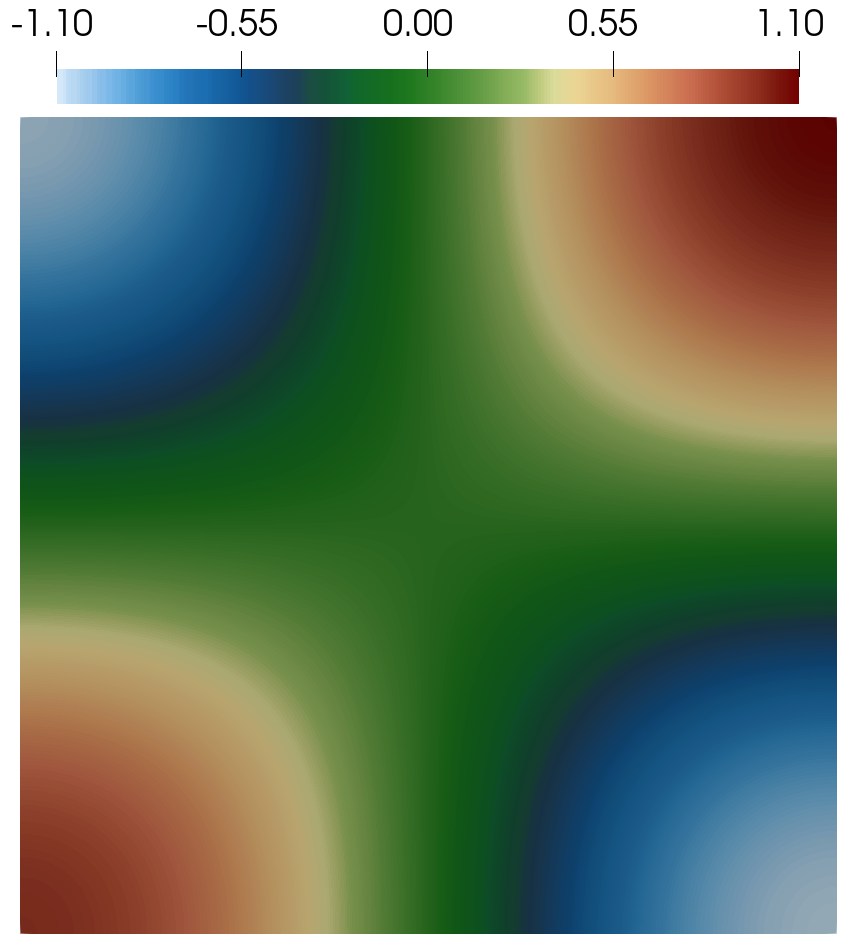} &
			\includegraphics[scale=.15]{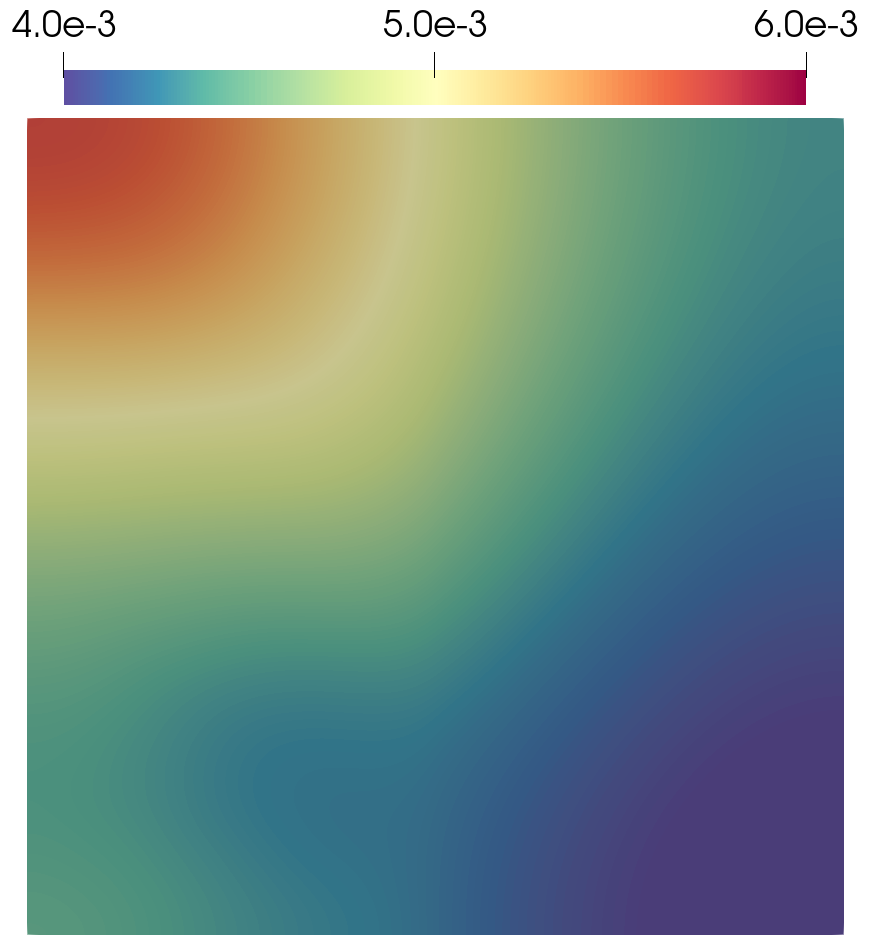}\\
			\includegraphics[scale=.15]{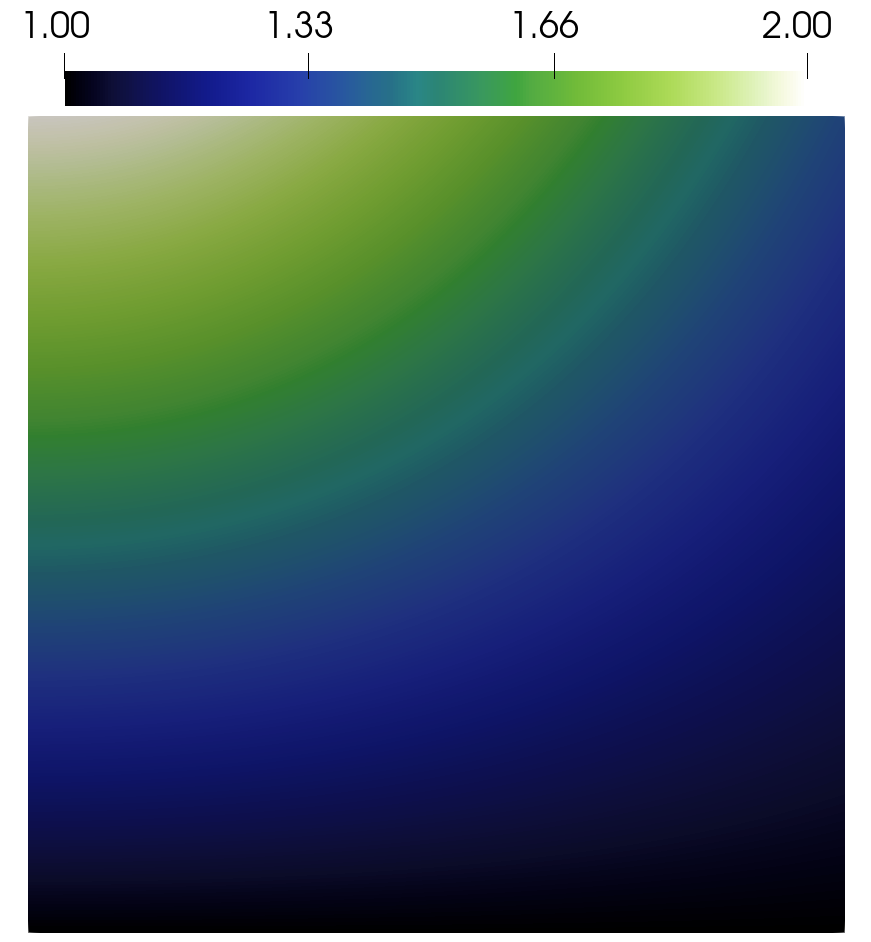} &
			\includegraphics[scale=.15]{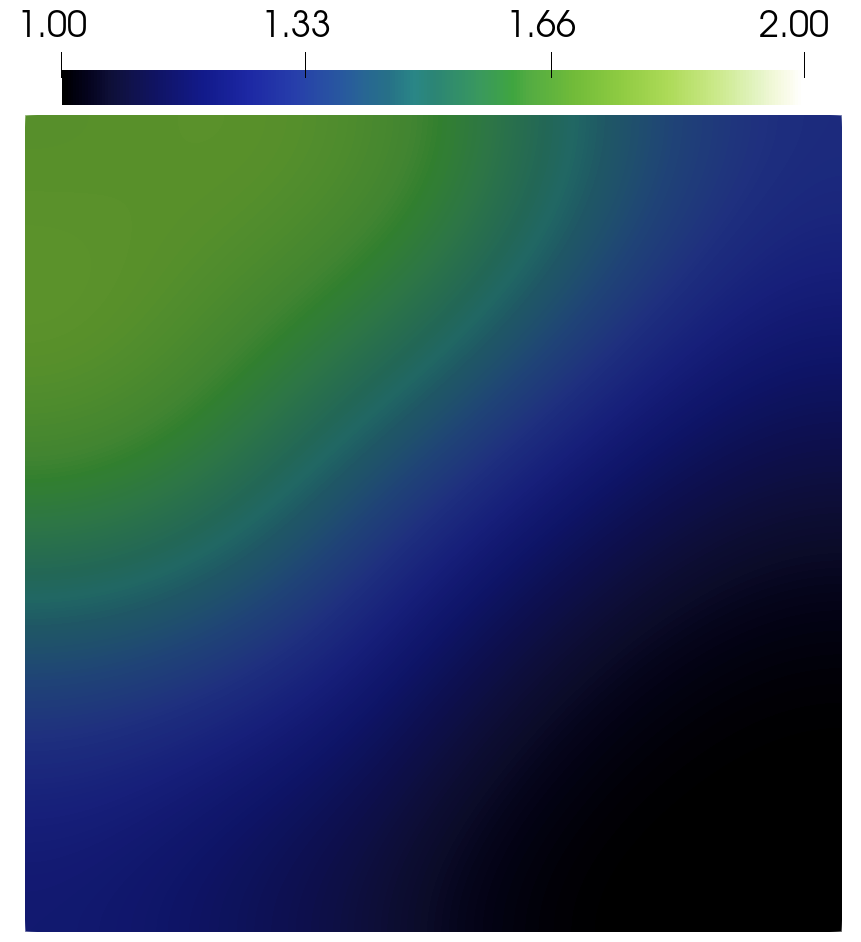} &
			\includegraphics[scale=.15]{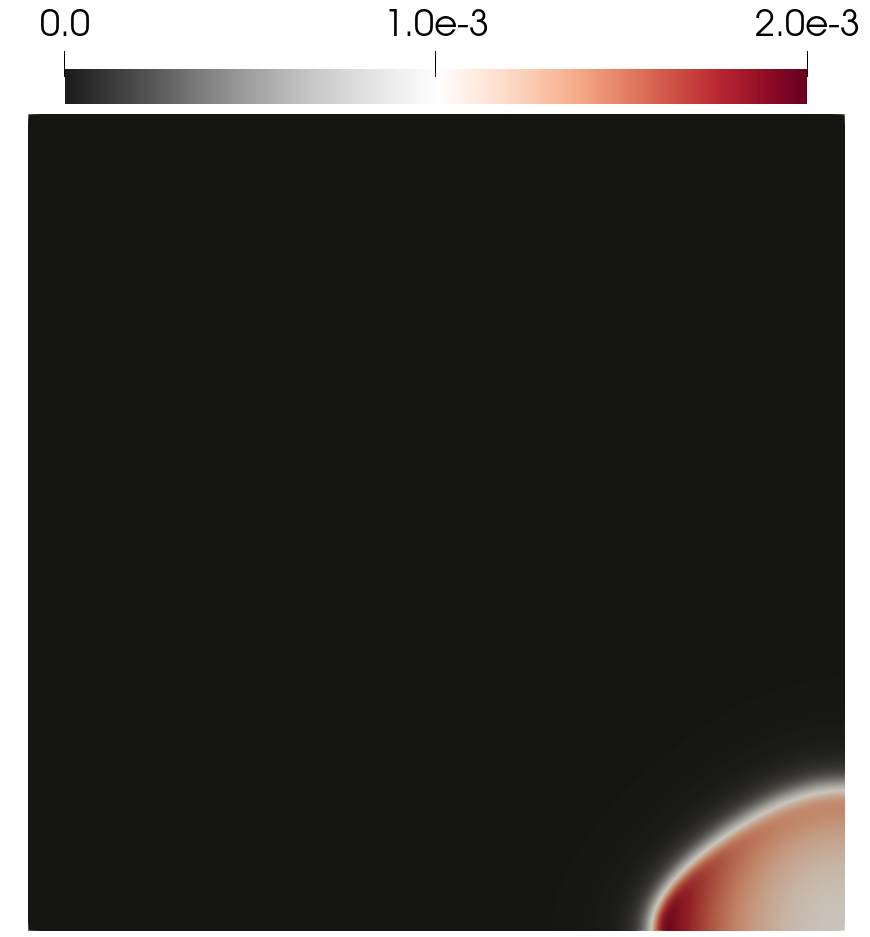} 
		\end{tabular}
		\caption{Top row: 
			noisy state observations $\statevar_{d,\zeta}$ (left),
			state reconstruction $\statevar^{\star}$ (middle),
			computed adjoint $\adjoint^{\star}$ (right).
			Bottom row: 
			true parameter $\param_{\text{true}}$ (left),
			parameter reconstruction $\param^{\star}$ (middle), 
			computed bound constraint Lagrange multiplier $\lboundmultiplier^{\star}$ (right). The dimension of each of the discretized fields is $124\,609$.}
		\label{fig:2Dreconstruction}
	\end{center}
\end{figure}

To provide greater context for the algorithmic and parallel scaling results presented later in this subsection we first present, in Figure~\ref{fig:2Dreconstruction}, the solution $(\statevar^{\star},\param^{\star})$,
also referred to as the reconstruction,
of a discretization of the example problem detailed in Section~\ref{sec:problemdescription}. Along with the state $\statevar^{\star}$, and parameter $\param^{\star}$, reconstructions, Figure~\ref{fig:2Dreconstruction} shows the noisy data $\statevar_{d,\zeta}$, the computed equality constraint Lagrange multiplier $\adjoint^{\star}$, and the computed bound constraint Lagrange multiplier $\lboundmultiplier^{\star}$. The significantly nonzero values of the bound constraint Lagrange multiplier indicate, when strict complementarity~\cite[Chapter 12]{nocedalwright2006} holds, at which points the lower-bound constraint $\param\geq\param_{\ell}$ is active i.e., where $\param^{\star}=\param_{\ell}$. 
In Table~\ref{tbl:algorithmicscaling}, the algorithmic scaling of the proposed framework is reported. By algorithmic scaling we mean asymptotic independence, with respect to mesh refinement, of the number of outer IP-Gauss-Newton steps to solve the optimization problem as well as the number of preconditioned Krylov subspace iterations required of each IP-Gauss-Newton linear system solve. 

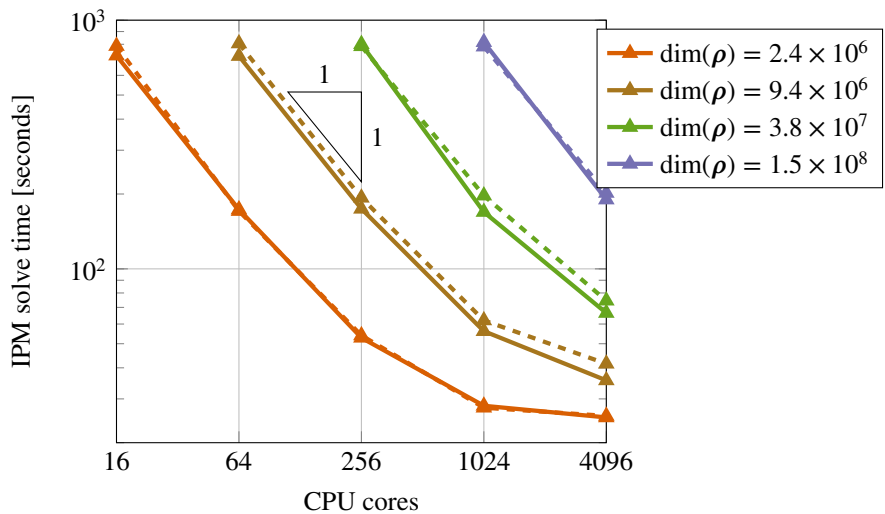
\begin{figure}[h!] 
	\begin{center}
		\begin{tikzpicture}
		\begin{axis}[
		grid=major, width=0.45\textwidth, height=0.4\textwidth,
		xlabel={CPU cores}, 
		ylabel={IPM solve time [seconds]},
		xmin=16, xmax=4096, ymin=2.e1, ymax=1.e3,
		xmode=log, ymode=log,
		xtick={16, 64, 256, 1024, 4096},
		xticklabels={16, 64, 256, 1024, 4096},
		ytick={1.e1, 1.e2, 1.e3, 1.e4},
		every axis plot/.append style={ultra thick},
		legend style={anchor=north west, nodes=right},
		xticklabel style={/pgf/number format/fixed},
		cycle list/Dark2
		]
		\addplot[index of colormap=1 of Dark2, mark=triangle*, mark options={solid}]
		table[x=cores, y=optTime] {GMRESscalingstrongScaling2.dat};
		\addlegendentry{dim$(\bparam)=2.4\times10^{6}$}
		\addplot[index of colormap=6 of Dark2, mark=triangle*, mark options={solid}]
		table[x=cores, y=optTime]      {GMRESscalingstrongScaling3.dat};
		\addlegendentry{dim$(\bparam)=9.4\times10^{6}$}
		\addplot[index of colormap=4 of Dark2, mark=triangle*, mark options={solid}]
		table[x=cores, y=optTime] {GMRESscalingstrongScaling4.dat};
		\addlegendentry{dim$(\bparam)=3.8\times10^{7}$}
		\addplot[index of colormap=2 of Dark2, mark=triangle*, mark options={solid}]
		table[x=cores, y=optTime] {GMRESscalingstrongScaling5.dat};
		\addlegendentry{dim$(\bparam)=1.5\times10^{8}$}
		\addplot[index of colormap=1 of Dark2, dashed, mark=triangle*, mark options={solid}]
		table[x=cores, y=optTime] {CGscalingstrongScaling2.dat};
		\addplot[index of colormap=6 of Dark2, dashed, mark=triangle*, mark options={solid}]
		table[x=cores, y=optTime]      {CGscalingstrongScaling3.dat};
		\addplot[index of colormap=4 of Dark2, dashed, mark=triangle*, mark options={solid}]
		table[x=cores, y=optTime]      {CGscalingstrongScaling4.dat};
		\addplot[index of colormap=2 of Dark2, dashed, mark=triangle*, mark options={solid}]
		table[x=cores, y=optTime]      {CGscalingstrongScaling5.dat};			
		\logLogSlopeTriangle{0.5}{0.15}{0.83}{-1}{black};
		\end{axis}
		\end{tikzpicture}
		\caption{
			Strong scaling for the IPM framework on the example problem~\ref{sec:problemdescription} using the preconditioned reduced-space IP-Gauss-Newton CG solver (dashed lines) described in Section~\ref{subsec:CGpreconditioner} and the Gauss-Seidel preconditioned IP-Gauss-Newton GMRES solver (solid lines) described in Section~\ref{subsec:GMRESpreconditioner}.}
		\label{fig:scaling}
	\end{center}
\end{figure}

Table~\ref{tbl:algorithmicscaling} provides numerical evidence, by the reported number of IP-Gauss-Newton linear system solves per optimizer computation, that the number of outer optimization steps is mesh independent. The average number of block Gauss-Seidel preconditioned GMRES iterations and $\bs{W}_{\! \bparam, \bparam}^{\mu}$ preconditioned CG iterations per IP-Gauss-Newton linear system solve shown in Table~\ref{tbl:algorithmicscaling} 
is numerical evidence that the performance of preconditioned Krylov subspace solvers as outlined in Sections~\ref{subsec:GMRESpreconditioner} and~\ref{subsec:CGpreconditioner} do not degrade under mesh refinement. 
\begin{table}[h!]
	\begin{center}
		\begingroup
		\setlength{\tabcolsep}{8pt} 
		\renewcommand{\arraystretch}{1.25} 
		\begin{tabular}{|c|c|c|c|}
			\hline
			& IP-Gauss-Newton 
			& Preconditioned GMRES & Preconditioned CG \\
			& linear solves per & iterations per & iterations per\\
			dim($\bparam$) & optimizer computation & IP-Gauss-Newton solve & 
			IP-Gauss-Newton solve
			\\ \hline
			$148\,609$      & $28.4$ & $6.50$ & $6.76$ 
			\\ \hline
			$591\,361$      & $28.2$ & $6.48$ & $6.72$ 
			\\ \hline
			$2\,362\,369$   & $28.8$ & $6.51$ & $6.68$
			\\ \hline
			$9\,443\,329$   & $28.3$ & $6.49$ & $6.85$
			\\ \hline
			$37\,761\,025$  & $28.7$ & $6.42$ & $6.87$
			\\ \hline
			$151\,019\,521$ & $29.0$ & $6.52$ & $6.75$
			\\ \hline 
		\end{tabular}
		\endgroup 
	\end{center} 
	\caption{Algorithmic scaling of the IP-Gauss-Newton method with block Gauss-Seidel preconditioned GMRES solves 
		and $\bs{W}_{\! \bparam, \bparam}^{\mu}$ preconditioned CG Schur complement solves of the IP-Gauss-Newton linear systems (see Equation~\eqref{eq:IPNewtonsys}). The absolute tolerance of the outer optimization loop is $10^{-6}$ and the relative tolerance of the block AMG-CG solves is $10^{-13}$.	
		CG terminates when $\|\bs{r}^{(k)}\|_{\bs{B}^{-1}} \leq \tau \|\bs{r}^{(0)}\|_{\bs{B}^{-1}}$ and GMRES terminates when $\|\bs{B}^{-1}\bs{r}^{(k)}\|_{2}\leq \tau \|\bs{B}^{-1}\bs{r}^{(0)}\|_{2}$, where $\tau=10^{-8}$ is the relative tolerance for the Krylov subspace solvers. Here, $\bs{B}$ denotes a generic preconditioner whose inverse approximates $\bs{A}^{-1}$.
		The degree of uniform mesh refinement is varied, as indicated by the dimension of the discretized parameter dim($\bparam$).}
	\label{tbl:algorithmicscaling}
\end{table}
Algorithmic scaling and the scalable performance of the AMG preconditioned CG solves for block matrices $\bs{J}_{\! \bstatevar}$, $\bs{J}_{\!\bstatevar}^{\top}$ and $\bs{W}_{\! \bparam, \bparam}^{\mu}$ that are involved in the application of the preconditioners, as detailed in Section~\ref{subsec:preconditionerCost} , provide for a PDE- and bound-constrained optimization method that has desirable parallel scaling, as reported in Figure~\ref{fig:scaling}.

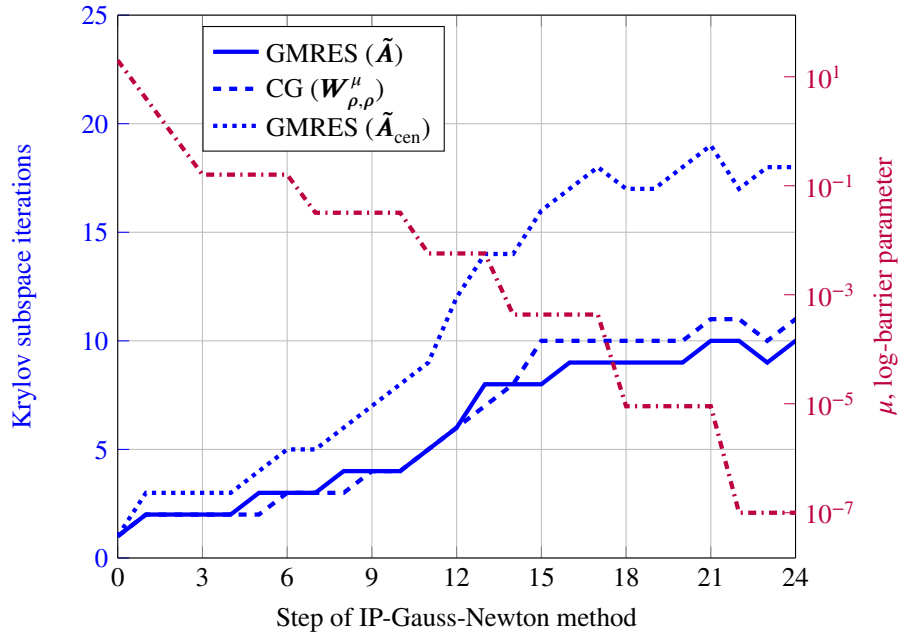
\begin{figure}[h!]
	\begin{center} 
		\begin{tikzpicture}
		\begin{axis}[
		scale only axis,
		width=0.5\textwidth, height=0.4\textwidth,
		xmin=0,xmax=24,
		ymin=0, ymax=25,
		grid=major,
		axis y line*=left,
		xlabel={Step of IP-Gauss-Newton method},
		ylabel={\textcolor{blue}{Krylov subspace iterations}},
		every axis plot/.append style={ultra thick},
		xtick={0,3,6,9,12,15,18,21,24},
		ytick={0, 5, 10, 15, 20, 25},
		yticklabels={\textcolor{blue}{$0$}, 
			\textcolor{blue}{$5$},
			\textcolor{blue}{$10$},
			\textcolor{blue}{$15$},
			\textcolor{blue}{$20$},
			\textcolor{blue}{$25$}
		},
		y tick style={color=blue},
		legend style={at={(0.125,0.98)}, anchor=north west, nodes=right},		
		]
		\addlegendentry{GMRES ($\bs{\tilde{A}}$)} 
		\addplot[blue, mark options={solid}]
		table[x=NewtonIt, y=KrylovIts] {BlockGSGMRESItsVsMu.dat};
		\addlegendentry{CG ($\bs{W}_{\! \bparam, \bparam}^{\mu}$)}
		\addplot[blue, dashed, mark options={solid}] table[x=NewtonIt, y=KrylovIts] {CGItsVsMu.dat};
		\addlegendentry{GMRES ($\bs{\tilde{A}}_{\text{cen}}$)} 
		\addplot[blue, dotted, mark options={solid}]
		table[x=NewtonIt, y=KrylovIts] {BlockJacobiGMRESItsVsMu.dat};
		\end{axis}
		\begin{axis}[
		scale only axis,
		width=0.5\textwidth, height=0.4\textwidth,
		xmin=0,xmax=24,
		axis y line*=right,
		axis x line=none,
		ymode=log,
		every axis plot/.append style={thick},
		ytick={1.e-7, 1.e-5, 1.e-3, 1.e-1, 1.e1},
		yticklabels={\textcolor{purple}{$10^{-7}$}, \textcolor{purple}{$10^{-5}$}, \textcolor{purple}{$10^{-3}$}, \textcolor{purple}{$10^{-1}$}, \textcolor{purple}{$10^{1}$}},
		y tick style={color=purple},
		ylabel={\textcolor{purple}{$\mu$, log-barrier parameter}},
		ylabel near ticks, yticklabel pos=right,
		]
		\addplot[purple, dashdotted, ultra thick, mark options={solid}]
		table[x=NewtonIt, y=Mu] {BlockJacobiGMRESItsVsMu.dat};
		\end{axis}
		\end{tikzpicture}
		\caption{
			The number of $\bs{W}_{\! \bparam, \bparam}^{\mu}$ preconditioned CG, block Gauss-Seidel ($\bs{\tilde{A}}$) preconditioned GMRES and central null ($\bs{\tilde{A}}_{\text{cen}}$) preconditioned GMRES iterations to solve the IP-Gauss-Newton linear system and the log-barrier $\mu$ at each step of the IP-Gauss-Newton method. The results are for the problem described in Section~\ref{sec:problemdescription} where the dimension of the discretized fields is equal to $2\,362\,369$. The number of Krylov subspace iterations has some initial variation, but quickly becomes relatively constant as the log-barrier parameter $\mu$ varies over approximately four orders of magnitude, where the central null preconditioner has a significantly higher iteration count. This provides numerical evidence that the preconditioned matrices have spectral properties that are asymptotically independent of $\mu$.}
		\label{fig:GMRESMINRESitsVsMu}
	\end{center}
\end{figure}

Figure~\ref{fig:GMRESMINRESitsVsMu} provides numerical evidence that supports Proposition~\ref{prop:GMRESresidualReduction},
that the performance of each proposed preconditioned Krylov subspace method is asymptotically independent of the log-barrier parameter $\mu$.
This feature of the preconditioned IP-Gauss-Newton Krylov subspace solvers is critical to enable efficient linear system solution estimation as IP-Newton linear systems are notoriously plagued by ill-conditioning due to the log-barrier parameter. In Figure~\ref{fig:GMRESMINRESitsVsMu} we compare the block Gauss-Seidel preconditioner against the central null preconditioner and see a significant improvement. Computational results presented in Table~\ref{tbl:algorithmicscaling}, Figure~\ref{fig:scaling}, and Figure~\ref{fig:GMRESMINRESitsVsMu} were generated on the Quartz cluster, with Intel Xeon E5-2695 v4 processors, at the Lawrence Livermore National Laboratory. The software extensively uses MFEM~\cite{andrej2024} version 4.5.1 and hypre~\cite{hypre} version 2.25.0.

\subsection{Regularization via the Morozov discrepancy principle} \label{subsec:Morozov}

\begin{figure}[h!] 
	\begin{center} 
		\begin{tabular}[c]{cc}
			\includegraphics[clip, scale=0.175]{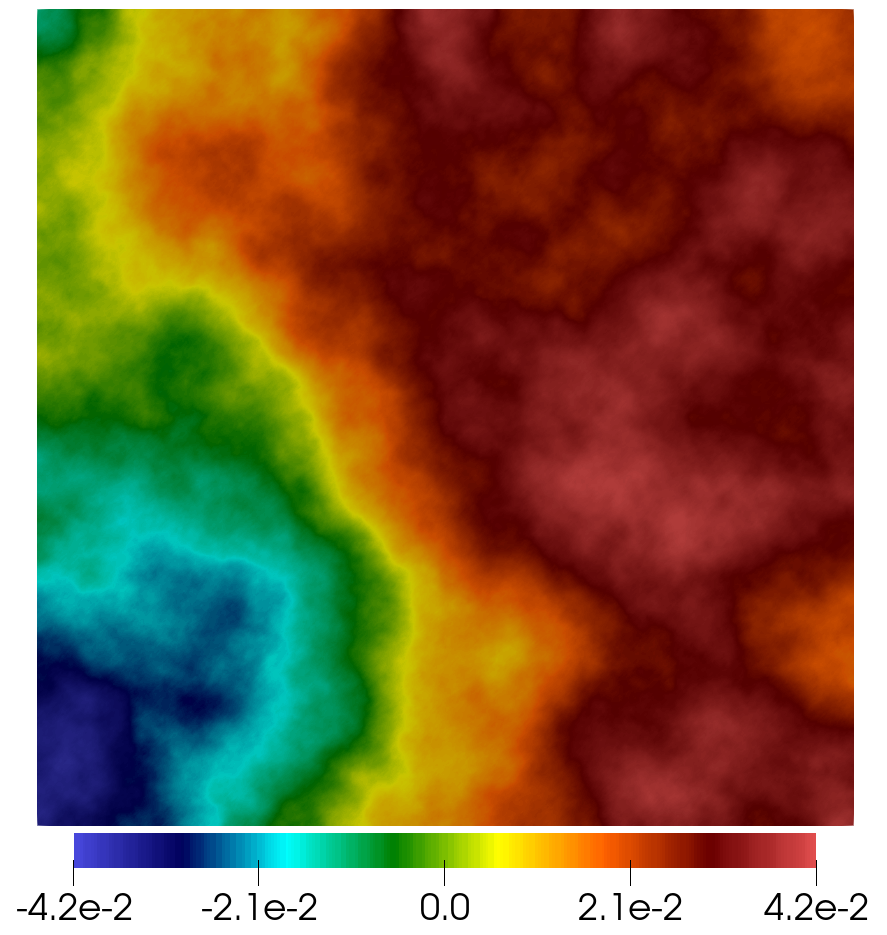}
			&
			\begin{tikzpicture}
			\begin{axis}[
			grid=major, width=0.5\textwidth, height=0.4\textwidth,
			xlabel={$\gamma_{1},\gamma_{2}$, regularization parameter}, 
			xmode=log, ymode=log,
			xmin=1.e-7, xmax=1.e-1,		
			xtick={1.e-7, 1.e-6, 1.e-5, 1.e-4, 1.e-3, 1.e-2, 1.e-1},
			every axis plot/.append style={ultra thick},
			legend style={at={(0.05, 0.95)}, anchor=north west, nodes=right},
			xticklabel style={/pgf/number format/fixed},
			cycle list/Dark2
			]
			\addplot[black!50!green, mark options={solid}]
			table[x=gamma, y=misfit] {MorozovData_noise0.05.dat};
			\addlegendentry{$\|\statevar^{\star}-\statevar_{d,\zeta}\|_{L^{2}(\Omega_{\text{left}})}$}
			\addplot[black, dashed, mark options={solid}]
			table[x=gamma, y=noisemisfit] {MorozovData_noise0.05.dat};
			\addlegendentry{$\|\zeta\|_{L^{2}(\Omega_{\text{left}})}$}
			\end{axis}
			\end{tikzpicture}
		\end{tabular}
	\end{center} 
	\caption{Left: spatial structure of a random sample $\zeta$. Right:
		Seminorm of discrepancy $(\statevar^{\star}-\statevar_{d,\zeta})$ and seminorm of noise $\zeta$ as functions of the regularization parameters $\gamma_{1}=\gamma_{2}$, plot obtained in order to apply the Morozov discrepancy principle for choosing numerical values of the regularization parameters with a $5\%$ noise level $\sigma_{\zeta}$.}
	\label{fig:Morozov}
\end{figure}
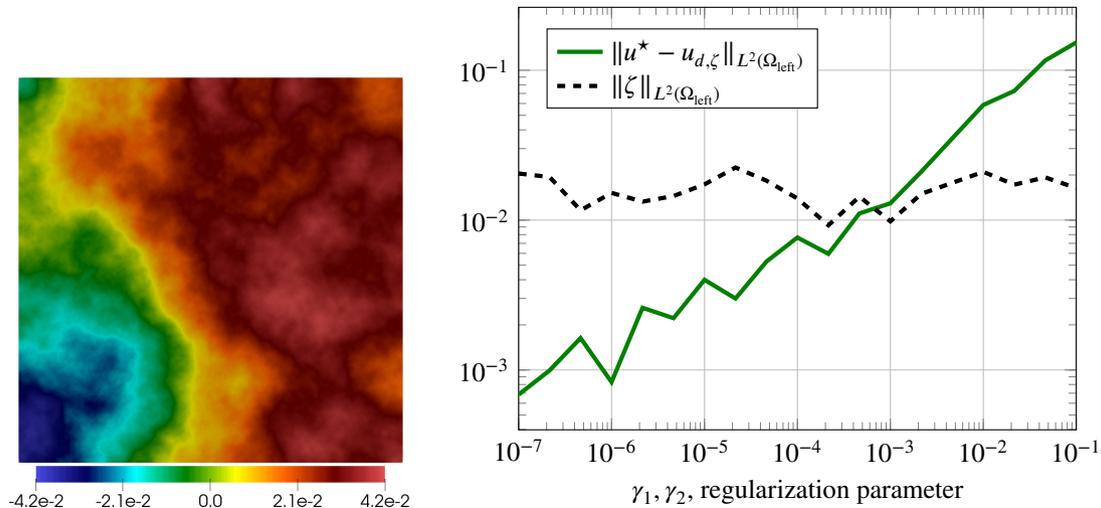

Choosing appropriate numerical values of the regularization parameters for ill-posed PDE-constrained optimization problems is essential, especially for those problems with noisy and incomplete data. If the regularization parameters, here $\gamma_{1},\gamma_{2}$, are too small then the optimal solution $(\statevar^{\star},\param^{\star})$ of~\eqref{eq:infdimproblem}, will be highly sensitive to the random noise $\zeta$. If the regularization parameters are too large, then the objective $\objfun$ does not place enough emphasis on $\objfun_{\text{misfit}}$ and $\statevar^{\star}$ will poorly resemble the data $\statevar_{d,\zeta}$. The Morozov discrepancy principle~\cite{morozov1966} is a means to choose the regularization parameters that avoids poorly chosen regularization parameters, thereby avoiding the aforementioned issues. Ultimately, for a given noise level, the regularization parameters are chosen so that the discrepancy $(\statevar^{\star}-\statevar_{d,\zeta})$, is equal to $\zeta$ in (semi)norm. For consistency with the data-misfit component of the objective $\objfun_{\text{misfit}}$ we use the seminorm $\|\cdot\|_{L^{2}(\Omega_{\text{left}})}$, i.e., the regularization parameters are chosen so that $\|\statevar^{\star}-\statevar_{d,\zeta}\|_{L^{2}(\Omega_{\text{left}})}=\|\zeta\|_{L^{2}(\Omega_{\text{left}})}$. For $5\%$ relative noise $\sigma_{\zeta}$, said equality is approximately achieved for $\gamma_{1}=\gamma_{2}=10^{-3}$ as seen in Figure~\ref{fig:Morozov} (right).       

In Table~\ref{tbl:noiseTable}, we report the average number of Gauss-Seidel preconditioned GMRES iterations required to solve the IP-Gauss-Newton linear systems for different noise levels $\sigma_{\zeta}$.
\begin{table}[h!]
	\begin{center}
		\begingroup
		\setlength{\tabcolsep}{6pt} 
		\renewcommand{\arraystretch}{1.25} 
		\begin{tabular}{|c|c|c|c|c|}
			\hline
			$\sigma_{\zeta}$ (noise level) & $1\%$ & $2\%$ & $5\%$ & $10\%$ \\ \hline 
			$\gamma_{1}$, $\gamma_{2}$ (Morozov) & 
			$2.2\times 10^{-4}$ & $4.6\times 10^{-4}$ &
			$1.0\times 10^{-3}$ & $2.2\times 10^{-3}$ \\ \hline 
			Mean Gauss-Seidel GMRES iterations & $7.9$ & $7.1$ & $6.6$ & $5.8$ \\ \hline 
		\end{tabular}
		\endgroup 
	\end{center}
	\caption{Regularization parameters $\gamma_{1}=\gamma_{2}$ as determined by the Morozov discrepancy principle for the noise level $\sigma_{\zeta}$ and the average number of block Gauss-Seidel preconditioned GMRES iterations required to solve the IP-Gauss-Newton systems to a~$10^{-8}$ relative error tolerance for a discretization such that dim$(\bparam)=2\,025$.} 
	\label{tbl:noiseTable}
\end{table}
Table~\ref{tbl:noiseTable} demonstrates
that the Gauss-Seidel preconditioner performs well for typical noise levels on the example problem detailed in Section~\ref{sec:problemdescription}. This performance is, however, not independent of the noise level. This dependent behavior is to be expected as, by the Morozov discrepancy principle, the regularization parameters increase with increasingly large noise levels. Larger regularization parameters in turn lead to a decrease in the eigenvalues of $\bs{H}_{\! \bparam, \bparam}^{-1}\bs{\hat{H}}$. By Proposition~\ref{prop:GMRESresidualReduction}, the number of block Gauss-Seidel preconditioned GMRES iterations required to solve IP-Gauss-Newton linear systems will decrease as the eigenvalues of $\bs{H}_{\! \bparam, \bparam}^{-1}\bs{\hat{H}}$ decrease, from above, to a value of one. In Table~\ref{tbl:regRobustnessTable}, we show the performance of the block Gauss-Seidel preconditioner over an increased range of regularization parameters, but with a fixed $5\%$ noise level.
\begin{table}[h!]
	\begin{center}
		\begingroup
		\setlength{\tabcolsep}{6pt} 
		\renewcommand{\arraystretch}{1.25} 
		\begin{tabular}{|c|c|c|c|c|c|}
			\hline
			$\gamma_{1}$, $\gamma_{2}$ & 
			$1.0\times 10^{-5}$ & $1.0\times 10^{-4}$ &
			$1.0\times 10^{-3}$ & $1.0\times 10^{-2}$ &
			$1.0\times10^{-1}$ \\ \hline 
			Mean Gauss-Seidel GMRES iterations & $12.4$ & $8.7$ & $6.6$ & $5.3$ & $4.7$ \\ \hline 
		\end{tabular}
		\endgroup 
	\end{center}
	\caption{The mean, with respect to each step of the IP-Gauss-Newton method, number of block Gauss-Seidel preconditioned GMRES iterations required to solve the IP-Gauss-Newton linear systems with a~$10^{-8}$ relative tolerance for the preconditioned residual with various regularization parameters $\gamma_{1}=\gamma_{2}$. All results in this table were generated with fixed $\sigma_{\zeta}=5\%$ and dim$(\bparam)=2\,025$.} 
	\label{tbl:regRobustnessTable}
\end{table}
%
%
As seen in Table~\ref{tbl:regRobustnessTable}, the block Gauss-Seidel preconditioner still performs well when the regularization parameters are chosen significantly smaller than they should.

\section{Conclusion}\label{sec:conclusion}
In this work, we presented a scalable means to solve a class of nonlinear elliptic PDE- and bound-constrained optimization problems described by Equation~\ref{eq:infdimproblem}. We applied a Newton-based interior-point method equipped with a globalizing filter line-search. The presented approach makes use of a Gauss-Newton search direction, for which the computational costs associated to inertia regularization, such as that from expensive additional linear solves, are avoided. We utilized a block Gauss-Seidel preconditioner and rigorously showed that for a class of PDE- and bound-constrained optimization problems that the eigenvalue clustering of the IP-Gauss-Newton matrix is asymptotically independent of discretization and ill-conditioning that arises from the IPM. We also presented a related preconditioned reduced-space Schur complement approach, that has identical eigenvalue clustering and for which CG is a valid Krylov subspace linear solver. For both of the aforementioned solvers, each preconditioned Krylov subspace iteration is made scalable by utilizing scalable AMG preconditioned CG solvers for subblock solves of the IP-Gauss-Newton matrix. We then demonstrated the scalable performance, with respect to mesh refinement and reduction of the log-barrier parameter, of the approach on a nonlinear elliptic PDE- and bound-constrained optimization problem.

\section{Acknowledgement}\label{sec:acknowledgement}

This work was performed under the auspices of the U.S. Department of Energy by Lawrence Livermore National Laboratory under Contract DE-AC52-07NA27344 and was supported by the LLNL-LDRD Program under Project No. 23-ERD-017. Release number LLNL-JRNL-858226-DRAFT. Support for this work was also provided by the National Science Foundation under Grant No. DMS-1840265 and CAREER-1654311.

\appendix

We next detail various supporting results, the first of which, Proposition~\ref{prop:eigenvalueOrdering} is a general linear algebraic result that we do not believe is entirely novel, but include it nonetheless for completeness as we make explicit use of it in Section~\ref{subsec:GMRESpreconditioner}. The second result, Proposition~\ref{prop:notDiagnolizable} shows that the block Gauss-Seidel preconditioned IP-Gauss-Newton matrix is not diagonalizable when the rank of $\bs{H}_{\!\bstatevar,\bstatevar}$ is less than the rank of $\bs{J}_{\!\bparam}$, which often is the case when $\bs{H}_{\!\bstatevar,\bstatevar}$ is rank-deficient. The third result, Proposition~\ref{prop:GSdiagonalization}, demonstrates explicitly that the perturbed Gauss-Seidel preconditioned IP-Gauss-Newton matrix is diagonalizable. 
Having established the diagonalizability of the perturbed preconditioned system, we then characterize the rate at which the residuals of the preconditioned GMRES solve in Proposition~\ref{prop:GMRESresidualReduction}.
In Appendix~\ref{sec:derivatives} we lastly list expressions for first and second order derivatives for the example problem described in Section~\ref{sec:problemdescription} as needed by the proposed method.

\section{Ordering of two sets of generalized eigenvalues}\label{sec:eigOrdering}

\begin{proposition} \label{prop:eigenvalueOrdering}
	Let $\bs{A}\in\mathbb{R}^{N\times N}$ be symmetric semidefinite and $\bs{B},\bs{C}\in\mathbb{R}^{N\times N}$ be symmetric positive definite where $(\bs{B}-\bs{C})$ is positive semidefinite. Consider the following generalized eigenvectors $\bs{u}^{(i)},\bs{v}^{(i)}$ and eigenvalues $\beta_{1}\geq\beta_{2}\dots\geq \beta_{N}\geq 0$, $\xi_{1}\geq \xi_{2}\geq\dots \geq \xi_{N}\geq 0$
	\begin{align*}
	\bs{A}\bs{u}^{(i)}&=\beta_{i}\bs{B}\bs{u}^{(i)}, \quad 1\leq i\leq N,\\
	\bs{A}\bs{v}^{(i)}&=\xi_{i}\bs{C}\bs{v}^{(i)}, \quad 1\leq i\leq N.
	\end{align*}
	Then, $\beta_{k}\leq \xi_{k}$ for each $k=1,2,\dots,N$.
	
	\begin{proof} 
		Let $k\in\lbrace 1,2,\dots,N\rbrace$ be arbitrary and define $E_{k}^{u}=$ Span$\lbrace \bs{u}^{(1)},\bs{u}^{(2)},\dots,\bs{u}^{(k)}\rbrace$ and $E_{k}^{v}=$ Span$\lbrace \bs{v}^{(1)},\bs{v}^{(2)},\dots,\bs{v}^{(k)}\rbrace$. 
		For $\bs{x}\in E_{k}^{u}$, there exists $\alpha_{1},\alpha_{2},\dots,\alpha_{k}\in\mathbb{R}$ such that $\bs{x}=\sum_{i=1}^{k}\alpha_{i}\bs{u}^{(i)}$. Then, utilizing the $\bs{B}$ orthogonality of $\lbrace \bs{u}^{(1)},\bs{u}^{(2)},\dots,\bs{u}^{(k)}\rbrace$ and that $\bs{B}\geq\bs{C}$ we have
		\begin{align*}
		\bs{x}^{\top}\bs{A}\bs{x}&=\sum_{i=1}^{k}\sum_{j=1}^{k}
		\alpha_{i}\alpha_{j}(\bs{u}^{(i)})^{\top}\bs{A}\bs{u}^{(j)}\\
		&=\sum_{i=1}^{k}\sum_{j=1}^{k}\alpha_{i}\alpha_{j}
		(\bs{u}^{(i)})^{\top}(\beta_{j}\bs{B}\bs{u}^{(j)})\\
		&=\sum_{i=1}^{k}\beta_{i}((\alpha_{i}\bs{u}^{(i)})^{\top}\bs{B}(\alpha_{i}\bs{u}^{(i)}))\\
		&\geq \beta_{k}\bs{x}^{\top}\bs{B}\bs{x}\\
		&\geq \beta_{k}\bs{x}^{\top}\bs{C}\bs{x}.
		\end{align*}
		We now define $E^{v,\perp}_{j}=\lbrace \bs{x}\in\mathbb{R}^{N}\setminus\lbrace\bs{0}\rbrace \text{ such that } \bs{x}^{\top}\bs{C}\bs{v}^{(i)}=0,\text{ for each }i\in\lbrace1,2,\dots,j\rbrace\rbrace$,
		$E^{u,\perp}_{j}=\lbrace \bs{x}\in\mathbb{R}^{N}\setminus\lbrace\bs{0}\rbrace \text{ such that } \bs{x}^{\top}\bs{B}\bs{u}^{(i)}=0,\text{ for each } i\in\lbrace1,2,\dots,j\rbrace\rbrace$. Utilizing these sets and the $\bs{C}$-orthogonality structure of the generalized eigenvectors $\bs{v}^{(1)},\bs{v}^{(2)},\dots\bs{v}^{(N)}$ we have
		\begin{align*}
		\xi_{k}=\sup_{\bs{x}\in E^{v,\,\perp}_{k-1}
		}\frac{\|\bs{x}\|_{\bs{A}}^{2}}{\|\bs{x}\|_{\bs{C}}^{2}}\geq \sup_{\bs{x}\in (E^{u}_{k}\cap E^{v,\,\perp}_{k-1})}\frac{\|\bs{x}\|_{\bs{A}}^{2}}{\|\bs{x}\|_{\bs{C}}^{2}}\geq \beta_{k}.
		\end{align*}
		We note that it is guaranteed that $(E^{u}_{k}\cap E^{v,\,\perp}_{k-1})\neq \emptyset$, as $E^{u}_{k}$ is a $k$-dimensional subspace of $\mathbb{R}^{N}$ and $E^{v,\,\perp}_{k-1}\cup\lbrace\bs{0}\rbrace$ is an $(N-(k-1))$-dimensional subspace of $\mathbb{R}^{N}$ and so $E^{u}_{k}\cap (E^{v,\,\perp}_{k-1}\cup\lbrace\bs{0}\rbrace)$ is a subspace of dimension at least one. Thus, the supremum over $(E^{u}_{k}\cap E^{v,\,\perp}_{k-1})$ is meaningful.
	\end{proof}
\end{proposition}

\section{Diagonalizability of the block Gauss-Seidel preconditioned IP-Gauss-Newton matrix} \label{sec:diagonalization}

Next, Proposition~\ref{prop:notDiagnolizable} shows that 
\begin{customnote}{I}
	\label{note:notDiagnolizable}
	the block Gauss-Seidel preconditioned IP-Gauss-Newton matrix $\bs{\tilde{A}}^{-1}\bs{A}$ 
	is not diagonalizable when $\text{Rank}(\bs{H}_{\! \bstatevar,\bstatevar})<\text{Rank}(\bs{J}_{\!\bparam})$,
\end{customnote}
\noindent which is the case for a large number of PDE-constrained optimization problems. Following that, it is shown in Proposition~\ref{prop:GSdiagonalization} that a small positive definite perturbation to $\bs{H}_{\!\bstatevar,\bstatevar}$ renders $\bs{\tilde{A}}^{-1}\bs{A}$ diagonalizable.

\begin{proposition}
	\label{prop:notDiagnolizable}
	Let $\bs{A}$ and $\bs{\tilde{A}}$ be specified by Equation~\ref{eq:IPNewtonsys} and Equation~\ref{eq:BlockGSPreconditioner} respectively,
	wherein $\text{Rank}(\bs{H}_{\! \bstatevar,\bstatevar})<\text{Rank}(\bs{J}_{\!\bparam})$, then the block Gauss-Seidel preconditioned IP-Gauss-Newton matrix $\bs{\tilde{A}}^{-1}\bs{A}$ is not diagonalizable.
	
	\begin{proof}
		First, consider
		\begin{align*}
		\bs{A}_{p}&=
		\begin{bmatrix}
		\bs{U} & \bs{V}^{\top} \\
		\bs{V} & \bs{W}
		\end{bmatrix},\,\,\,
		\bs{U}=
		\begin{bmatrix}
		\bs{H}_{\! \bstatevar, \bstatevar} & \bs{J}_{\! \bstatevar}^{\top} \\
		\bs{J}_{\! \bstatevar} & \bs{0}
		\end{bmatrix}, \\
		\bs{V}&=\begin{bmatrix}
		\bs{0} & \bs{J}_{\! \bparam}^{\top} 
		\end{bmatrix},\,\,\,\,
		\bs{W}=\bs{W}_{\! \bparam, \bparam}^{\mu},
		\end{align*}
		where $\bs{A}_{p}$ is obtained from $\bs{A}$ by a symmetric permutation of its rows and columns. When expressed in this basis the block Gauss-Seidel preconditioner is block lower triangular
		\begin{align} \label{eq:BlockGSupperTriangular}
		\bs{\tilde{A}}_{p}&=
		\begin{bmatrix}
		\bs{U} & \bs{0} \\
		\bs{V} & \bs{W}
		\end{bmatrix}.
		\end{align}
		The block Gauss-Seidel preconditioned IP-Gauss-Newton matrix is then
		\begin{align*}
		\bs{\tilde{A}}_{p}^{-1}\bs{A}_{p}&=
		\begin{bmatrix}
		\bs{I} & \bs{U}^{-1}\bs{V}^{\top} \\
		\bs{0} & \bs{I}-\bs{W}^{-1}\bs{V}\bs{U}^{-1}\bs{V}^{\top}
		\end{bmatrix}.
		\end{align*}
		
		Given the block triangular structure of $\bs{\tilde{A}}_{p}^{-1}\bs{A}_{p}$ and that $\bs{I}-\bs{W}^{-1}\bs{V}\bs{U}^{-1}\bs{V}=\bs{I}+\bs{W}^{-1}\bs{\hat{H}}_{\bs{d}}$,  we see that $\bs{\tilde{A}}_{p}^{-1}\bs{A}_{p}$ contains $2\dimstate+\dimparam-\text{Rank}(\bs{\hat{H}}_{\!\bs{d}})$ eigenvalues that are equal to $1$. 
		We next show that the geometric multiplicity associated to this eigenvalue is less than the algebraic multiplicity. Towards this end we enumerate the linearly independent vectors that satisfy 
		\begin{align}
		\label{eq:eig1}
		\begin{bmatrix}
		\bs{I} & \bs{U}^{-1}\bs{V}^{\top} \\
		\bs{0} & \bs{I}-\bs{W}^{-1}\bs{V}\bs{U}^{-1}\bs{V}^{\top}
		\end{bmatrix}
		\begin{bmatrix}
		\bs{v}_{1} \\
		\bs{v}_{2}
		\end{bmatrix}
		&=
		1\begin{bmatrix}
		\bs{v}_{1} \\
		\bs{v}_{2}
		\end{bmatrix}.
		\end{align}
		There are $2\dimstate$ independent eigenvectors for which $\bs{v}_{2}=\bs{0}$, namely the standard unit vectors $\bs{v}_{1}=\bs{e}_{1},\bs{e}_{2},\dots,\bs{e}_{2\dimstate}$. If $\bs{\tilde{A}}_{p}^{-1}\bs{A}_{p}$ is nondefective then there must be $\dimparam-\text{Rank}(\bs{\hat{H}}_{\!\bs{d}})=\text{dim}(\text{Ker}(\bs{\hat{H}}_{\!\bs{d}}))$ linearly independent eigenvectors, with eigenvalue $\lambda=1$, for which $\bs{v}_{2}\neq\bs{0}$. 
		Since $\bs{U}$ is invertible, the first block row of Equation~\eqref{eq:eig1} then requires that $\bs{v}_{2}$ is a nonzero element of $\text{Ker}(\bs{V}^{\top})$. Given the block structure of $\bs{V}$ we can further conclude that $\bs{v}_{2}$ is a nonzero element of $\text{Ker}(\bs{J}_{\!\bparam})$. However, there are not $\text{dim}(\text{Ker}(\bs{\hat{H}}_{\!\bs{d}}))$  linearly independent vectors 
		as $\text{dim}(\text{Ker}(\bs{\hat{H}}_{\!\bs{d}}))\geq \text{dim}(\text{Ker}(\bs{H}_{\!\bstatevar,\bstatevar}))$ and by the assumptions of the Proposition, $\text{dim}(\text{Ker}((\bs{H}_{\! \bstatevar,\bstatevar}))>\text{dim}(\text{Ker}(\bs{J}_{\!\bparam}))$.			  
	\end{proof}	
\end{proposition}

\begin{proposition} \label{prop:GSdiagonalization}
	Let $\bs{A}_{\varepsilon}$ and $\bs{\tilde{A}}_{\varepsilon}$ be specified by
	\begin{align}
	\label{eq:perturbedA}
	\bs{A}_{\varepsilon}
	=
	\begin{bmatrix}
	\bs{H}_{\!\bstatevar,\bstatevar}+\varepsilon\bs{M}_{\!\bstatevar} & \bs{0} & \bs{J}_{\!\bstatevar}^{\top} \\
	\bs{0} & \bs{W}_{\!\bparam,\bparam}^{\mu} & \bs{J}_{\!\bparam}^{\top} \\
	\bs{J}_{\!\bstatevar} & \bs{J}_{\!\bparam} & \bs{0}
	\end{bmatrix},\,\,\,	\bs{\tilde{A}}_{\varepsilon}
	=
	\begin{bmatrix}
	\bs{H}_{\!\bstatevar,\bstatevar}+\varepsilon\bs{M}_{\!\bstatevar} & \bs{0} & \bs{J}_{\!\bstatevar}^{\top} \\
	\bs{0} & \bs{W}_{\!\bparam,\bparam}^{\mu} & \bs{J}_{\!\bparam}^{\top} \\
	\bs{J}_{\!\bstatevar} & \bs{0} & \bs{0}
	\end{bmatrix},	
	\end{align}
	where the scalar $\varepsilon>0$ is arbitrary, the mass matrix $\bs{M}_{\!\bstatevar}$ is positive definite, and $\bs{J}_{\!\bparam}$ is assumed to be full rank, then the perturbed block Gauss-Seidel preconditioned IP-Gauss-Newton matrix $\bs{\tilde{A}}_{\varepsilon}^{-1}\bs{A}_{\varepsilon}$ is diagonalizable.
	
	\begin{proof}
		First, as in the proof of Proposition~\ref{prop:notDiagnolizable}, consider
		\begin{align*}
		\bs{A}_{\varepsilon,p}&=
		\begin{bmatrix}
		\bs{U}_{\! \varepsilon} & \bs{V}^{\top} \\
		\bs{V} & \bs{W}
		\end{bmatrix},\,\,\,
		\bs{U}_{\! \varepsilon}=
		\begin{bmatrix}
		\bs{H}_{\! \bstatevar, \bstatevar}+\varepsilon\bs{M}_{\!\bstatevar} & \bs{J}_{\! \bstatevar}^{\top} \\
		\bs{J}_{\! \bstatevar} & \bs{0}
		\end{bmatrix}, \\
		\bs{V}&=\begin{bmatrix}
		\bs{0} & \bs{J}_{\! \bparam}^{\top} 
		\end{bmatrix},\,\,\,\,
		\bs{W}=\bs{W}_{\! \bparam, \bparam}^{\mu},
		\end{align*}
		where $\bs{A}_{\varepsilon,p}$ is obtained from $\bs{A}_{\varepsilon}$ by a symmetric permutation of its rows and columns. When expressed in this basis the perturbed block Gauss-Seidel preconditioner is block lower triangular
		\begin{align} 
		\bs{\tilde{A}}_{\varepsilon,p}&=
		\begin{bmatrix}
		\bs{U}_{\! \varepsilon} & \bs{0} \\
		\bs{V} & \bs{W}
		\end{bmatrix}.
		\end{align}
		The perturbed block Gauss-Seidel preconditioned IP-Gauss-Newton matrix is then
		\begin{align*}
		\bs{\tilde{A}}_{\varepsilon,p}^{-1}\bs{A}_{\varepsilon,p}&=
		\begin{bmatrix}
		\bs{I} & \bs{U}_{\! \varepsilon}^{-1}\bs{V}^{\top} \\
		\bs{0} & \bs{I}-\bs{W}^{-1}\bs{V}\bs{U}_{\! \varepsilon}^{-1}\bs{V}^{\top}
		\end{bmatrix}.
		\end{align*}
		Let $\bs{Q}\bs{\Lambda}\bs{Q}^{\top}=-\bs{W}^{-1/2}\bs{V}\bs{U}_{\! \varepsilon}^{-1}\bs{V}^{\top}\bs{W}^{-1/2}$ be the (symmetric) Schur decomposition~\cite[Theorems 7.1.3., 8.1.1.]{golub2013}, where $\bs{Q}$ is orthogonal $(\bs{Q}^{-1}=\bs{Q}^{\top})$ and $\bs{\Lambda}$ is diagonal. Furthermore, let $\bs{X}=\bs{U}_{\! \varepsilon}^{-1}\bs{V}^{\top}\bs{W}^{-1/2}\bs{Q}\bs{\Lambda}^{-1}$ then
		\begin{align*}
		\begin{bmatrix}
		\bs{I} & \bs{X} \\
		\bs{0} & \bs{W}^{-1/2}\bs{Q}
		\end{bmatrix}
		\begin{bmatrix}
		\bs{I} & \bs{0} \\
		\bs{0} & \bs{I}+\bs{\Lambda} 
		\end{bmatrix}
		\begin{bmatrix}
		\bs{I} & \bs{X} \\
		\bs{0} & \bs{W}^{-1/2}\bs{Q}
		\end{bmatrix}^{-1}&=\\
		\begin{bmatrix}
		\bs{I} & \bs{X} \\
		\bs{0} & \bs{W}^{-1/2}\bs{Q}
		\end{bmatrix}
		\begin{bmatrix}
		\bs{I} & \bs{0} \\
		\bs{0} & \bs{I}+\bs{\Lambda} 
		\end{bmatrix}
		\begin{bmatrix}
		\bs{I} & -\bs{X}\bs{Q}^{\top}\bs{W}^{1/2} \\
		\bs{0} & \bs{Q}^{\top}\bs{W}^{1/2}
		\end{bmatrix}&=\\
		\begin{bmatrix}
		\bs{I} & \bs{X}(\bs{I}+\bs{\Lambda}) \\
		\bs{0} & \bs{W}^{-1/2}\bs{Q}(\bs{I}+\bs{\Lambda})
		\end{bmatrix}
		\begin{bmatrix}
		\bs{I} & -\bs{X}\bs{Q}^{\top}\bs{W}^{1/2} \\
		\bs{0} & \bs{Q}^{\top}\bs{W}^{1/2}
		\end{bmatrix}&=\\
		\begin{bmatrix}
		\bs{I} & \bs{X}\bs{\Lambda}\bs{Q}^{\top}\bs{W}^{1/2}\\
		\bs{0} & \bs{I} + \bs{W}^{-1/2}\bs{Q}\bs{\Lambda}\bs{Q}^{\top}\bs{W}^{1/2}
		\end{bmatrix}&=\\
		\begin{bmatrix}
		\bs{I} & \bs{U}_{\! \varepsilon}^{-1}\bs{V}^{\top}
		\\
		\bs{0} & \bs{I} -\bs{W}^{-1}\bs{V}\bs{U}_{\! \varepsilon}^{-1}\bs{V}^{\top}
		\end{bmatrix}&=\\
		\bs{\tilde{A}}_{\varepsilon,p}^{-1}\bs{A}_{\varepsilon,p}
		\end{align*}
	\end{proof}
\end{proposition}
We note that each column of
\begin{align} \label{eq:GSeigVectors}
\bs{Y}&=\begin{bmatrix}
\bs{I} & \bs{X} \\
\bs{0} & \bs{W}^{-1/2}\bs{Q}
\end{bmatrix},
\end{align}
is an eigenvector of $\bs{\tilde{A}}_{\varepsilon,p}^{-1}\bs{A}_{\varepsilon,p}$, which is particularly relevant for analyzing the rate of convergence of a perturbed block Gauss-Seidel preconditioned GMRES solve of a perturbed IP-Gauss-Newton linear system as in Appendix~\ref{sec:perturbedIPGaussNewton}.

\section{Convergence rate of perturbed block Gauss-Seidel IP-Gauss-Newton GMRES solves}
\label{sec:perturbedIPGaussNewton}

In this section, the reduction in preconditioned residual norms associated to a perturbed block Gauss-Seidel preconditioned GMRES solve of a perturbed IP-Gauss-Newton linear system is characterized.

\begin{proposition} \label{prop:GMRESresidualReduction}
	The perturbed IP-Gauss-Newton linear system residuals $\bs{r}^{(k)}=\bs{b}-\bs{A}_{\varepsilon}\bs{x}^{(k)}$, generated from a perturbed block Gauss-Seidel $\bs{\tilde{A}}_{\varepsilon}$ (see Equation~\eqref{eq:perturbedA}) preconditioned GMRES solve as they arise in a Gauss-Newton-IPM for an elliptic PDE- and bound-constrained optimization problem, with a positive semidefinite subblock $\bs{H}_{\! \bstatevar, \bstatevar}$, symmetric positive definite $\bs{H}_{\! \bparam, \bparam},\bs{M}_{\!\bstatevar}$ and positive $\varepsilon$ satisfy
	\begin{align} \label{eq:residualReduction}
	\frac{\|\bs{\tilde{A}}_{\varepsilon}^{-1}\bs{r}^{(k+1)}\|_{2}}{
		\|\bs{\tilde{A}}_{\varepsilon}^{-1}\bs{r}^{(0)}\|_{2}}		
	\leq \delta_{k}\,\kappa(\bs{Y}),
	\end{align}   
	where $\bs{Y}$ is a permutation of the matrix of right eigenvectors of $\bs{\tilde{A}}_{\varepsilon}^{-1}\bs{A}_{\varepsilon}$ as specified in Equation~\eqref{eq:GSeigVectors}. The scalar $\delta_{k}$, is asymptotically independent of both discretization and log-barrier parameter $\mu$, and is given by
	\begin{align*}
	\delta_{k}=		\prod_{j=1}^{k}\left(\frac{\lambda_{j}(\bs{H}_{\! \bparam, \bparam}^{-1}\bs{\hat{H}}_{\bs{d},\varepsilon})}{1+\lambda_{j}(\bs{H}_{\! \bparam, \bparam}^{-1}\bs{\hat{H}}_{\bs{d},\varepsilon})}\right),
	\end{align*}
	where $\bs{\hat{H}}_{\bs{d},\varepsilon}=(\bs{J}_{\!\bstatevar}^{-1}\bs{J}_{\!\bparam})^{\top}(\bs{H}_{\!\bstatevar,\bstatevar}+\varepsilon\bs{M}_{\!\bstatevar})(\bs{J}_{\!\bstatevar}^{-1}\bs{J}_{\!\bparam})$
	For those problems that satisfy Property~\ref{propertyP}, $\delta_{k}$ is independent of discretization and ill-conditioning from the IPM.

	\begin{proof}
		~\cite[Proposition 4]{saad1986} along with the diagonalizability of $\bs{\tilde{A}}_{\varepsilon}^{-1}\bs{A}_{\varepsilon}$, as demonstrated in Proposition~\ref{prop:GSdiagonalization} of Appendix~\ref{sec:diagonalization}, allows us to bound the norm of the preconditioned residual $\|\bs{\tilde{A}}_{\varepsilon}^{-1}\bs{r}^{(k+1)}\|_{2}$ at the $k$th step of a GMRES solve of a IP-Gauss-Newton linear system as
		\begin{align} \label{eq:residualReductionSaad}
		\frac{\|\bs{\tilde{A}}_{\varepsilon}^{-1}\bs{r}^{(k+1)}\|_{2}}{
			\|\bs{\tilde{A}}_{\varepsilon}^{-1}\bs{r}^{(0)}\|_{2}}
		\leq \kappa(\bs{Y})\min_{p\in\mathbb{P}_{k},p(0)=1}\max_{1\leq j\leq \dimparam+2\dimstate}p(\lambda_{j}(\bs{\tilde{A}}_{\varepsilon}^{-1}\bs{A}_{\varepsilon})),
		\end{align}   
		Next, as in~\cite{saad1986}, we take $q\in\mathbb{P}_{k}$ given by
		\begin{align*}
		q(\lambda)=\prod_{j=1}^{k}q_{j}(\lambda),\,\,\,q_{j}(\lambda)= \left(1-\frac{\lambda}{\lambda_{j}}\right),
		\end{align*}
		where for notational simplicity we have used the symbol $\lambda_{j}=\lambda_{j}((\bs{W}_{\! \bparam, \bparam}^{\mu})^{-1}\bs{\hat{H}}_{\varepsilon})=1+\lambda_{j}((\bs{W}_{\! \bparam, \bparam}^{\mu})^{-1}\bs{\hat{H}}_{\bs{d},\varepsilon})$ and will continue to do so for the remainder of this proof. We then note that
		\begin{align*}
		\min_{p\in\mathbb{P}_{k},p(0)=1}\max_{1\leq j\leq \dimparam+2\dimstate}p(\lambda_{j})\leq \max_{1\leq j\leq \dimparam+2\dimstate}q(\lambda_{j})=
		\max_{k+1\leq j\leq \dimparam+2\dimstate}q(\lambda_{j}),
		\end{align*}
		and that $q_{j}(\lambda)$ is a nonnegative decreasing function of $\lambda$ on $[1,\lambda_{j}]\subset[1, \max_{k+1\leq i\leq \dimparam+2\dimstate}\lambda_{i}]$, so that
		\begin{align*}
		\max_{k+1\leq j\leq \dimparam+2\dimstate}q(\lambda_{j})\leq q(1)&=\prod_{j=1}^{k}\left(1-\frac{1}{\lambda_{j}}\right)=\prod_{j=1}^{k}\left(\frac{\lambda_{j}-1}{\lambda_{j}}\right)
		\\
		&=\prod_{j=1}^{k}\left(\frac{\lambda_{j}((\bs{W}_{\! \bparam, \bparam}^{\mu})^{-1}\bs{\hat{H}}_{\bs{d},\varepsilon})}{1+\lambda_{j}((\bs{W}_{\! \bparam, \bparam}^{\mu})^{-1}\bs{\hat{H}}_{\bs{d},\varepsilon})}\right) \leq 
		\prod_{j=1}^{k}\left(\frac{\lambda_{j}(\bs{H}_{\! \bparam, \bparam}^{-1}\bs{\hat{H}}_{\bs{d},\varepsilon})}{1+\lambda_{j}(\bs{H}_{\! \bparam, \bparam}^{-1}\bs{\hat{H}}_{\bs{d},\varepsilon})}\right)=\delta_{k}.
		\end{align*}
		The final inequality is due to $\lambda/(1+\lambda)$ being a nonnegative increasing function of $\lambda$ on $\mathbb{R}_{\geq 0}$. 
		To characterize the upper-bound $\delta_{k}$, we recall that Property~\ref{propertyP} holds for a large class of PDE-constrained optimization problems and for which the eigenvalues of 
		$\bs{H}_{\! \bparam, \bparam}^{-1}\bs{\hat{H}}_{\bs{d},\varepsilon}$ decay rapidly to zero and in a discretization independent manner, and since $\delta_{k}$ is defined only in terms of the eigenvalues of $\bs{H}_{\! \bparam, \bparam}^{-1}\bs{\hat{H}}_{\bs{d},\varepsilon}$ we can conclude that it also is asymptotically independent of discretization and log-barrier parameter~$\mu$. 		
	\end{proof}
\end{proposition}

\begin{customremark}{1}
	\label{rmk:perturbation}
	The matrix $\bs{\hat{H}}_{\bs{d},\varepsilon}$ obtained by perturbing $\bs{H}_{\!\bstatevar,\bstatevar}$, is the reduced-space IP-Gauss-Newton Hessian that one would obtain by adding $\varepsilon/2\int_{\Omega}(u-u_{d,\zeta})^{2}\mathrm{d}V$ to the objective functional that defines the PDE- and bound-constrained optimization problem in Equation~\eqref{eq:infdimproblem}. This small perturbation does not impact the qualitative features of $\bs{\hat{H}}_{\bs{d}}$ and so Property~\ref{propertyP} still holds. 
\end{customremark}

\begin{customremark}{2}
	\label{rmk:perturbedresiduals}
	Proposition~\ref{prop:GMRESresidualReduction} allows us to conclude that for $k$ large-enough, but independent of both the discretization and log-barrier parameter $\mu$, that $\delta_{k}\ll 1$. If we then assume that, $\kappa(\bs{Y})$ has, at worst, mild dependence on the mesh and log-barrier parameter, then by Equation~\eqref{eq:residualReduction}, a small and largely mesh- and log-barrier parameter-independent number of perturbed block Gauss-Seidel preconditioned GMRES iterations are required to achieve a specified relative reduction of the preconditioned residual-norm for the perturbed IP-Gauss-Newton linear system. To improve the proposed preconditioner, one could use a better approximation of the Schur-complement $\bs{\hat{H}}=\bs{\hat{H}}_{\bs{d}}+\bs{W}_{\! \bparam, \bparam}^{\mu}$ than $\bs{W}_{\! \bparam, \bparam}^{\mu}$, which would require a means to approximate $\bs{\hat{H}}_{\bs{d}}$. However, it is challenging to approximate $\bs{\hat{H}}_{\bs{d}}$, since it is a formally dense matrix free operator defined (see Equation~\eqref{eq:reduced_misfithessian}) in terms of inverse elliptic PDE operators, e.g., $\bs{J}_{\! \bstatevar}^{-1}$. Since, $\bs{\hat{H}}_{\bs{d}}$ has structure such as rapidly decaying eigenvalues that allow for global low-rank approximations~\cite{buighattasetal2013}. Due to narrow sensitivities of the parameter-to-observable map~\cite[Section 3.1]{hartland2023} $\bs{\hat{H}}_{\bs{d}}$ also has hierarchical low-rank structure~\cite{ambartsumyan2020, alger2024point}. Quasi-Newton approximations are also available, see e.g.,~\cite{birosghattas2005, vuchkov2020}.
	We leave these potential improvements for future work.
\end{customremark}

\section{Derivatives for the nonlinear elliptic PDE- and bound-constrained example problem}\label{sec:derivatives}
For completeness we now list various derivatives of the various discretized quantites that are contained in Section~\ref{subsec:interiorpoint} for the example problem described in Section~\ref{sec:problemdescription}.  First, we recall that as in Section~\ref{subsec:interiorpoint} $\statevar,\param$ are approximated by elements from the finite-dimensional spaces $\mathcal{V}_{h},\mathcal{M}_{h}$ respectively with basis sets $\lbrace\phi_{j}(\spatialvar)\rbrace_{j=1}^{\dimstate}\subset\mathcal{V}_{h}$ and $\lbrace \psi_{j}(\spatialvar)\rbrace_{j=1}^{\dimparam}\subset\mathcal{M}_{h}$. We begin with the discretized PDE-constraint $\bs{c}(\bstatevar,\bparam)$ described in Equation~\eqref{eq:examplePDEconstraint} and the Jacobians $\bs{J}_{\! \bstatevar},\bs{J}_{\! \bparam}$ of it with respect to the state $\bstatevar$, and parameter $\bparam$,
\begin{align*} 
(\bs{c}(\bstatevar,\bparam))_{i}&=\int_{\Omega}\left(\rho\bs{\nabla}_{\! \spatialvar}\phi_{i}\cdot\bs{\nabla}_{\! \spatialvar}\statevar 
+\phi_{i}(\statevar+\statevar^{3}/3-g)\right)\mathrm{d}\spatialvar,\,\,\,1\leq i\leq \dimstate, \\
(\bs{J}_{\! \bstatevar})_{i,j}&=\int_{\Omega}\left(\rho\bs{\nabla}_{\! \spatialvar}\phi_{i}\cdot\bs{\nabla}_{\! \spatialvar}\phi_{j} 
+\phi_{i}\phi_{j}(1+\statevar^{2})\right)\mathrm{d}\spatialvar,\,\,\,1\leq i,j\leq \dimstate,\\
(\bs{J}_{\! \bparam})_{i,j}&=\int_{\Omega}\left(\psi_{j}\bs{\nabla}_{\! \spatialvar}\statevar\cdot\bs{\nabla}_{\! \spatialvar}\phi_{i}\right) \mathrm{d}\spatialvar,\,\,\,1\leq i\leq \dimstate,\,\,\,1\leq j\leq \dimparam.\\
\end{align*}
Finally, we list gradients and Hessians of the objective $\objfun_{h}$ that are relevant for a Gauss-Newton approach
\begin{align*}
(\bs{\nabla}_{\! \bstatevar}\objfun_{h})_{i}&=
\int_{\Omega_{\text{left}}}\phi_{i}(\statevar_{h}-\statevar_{d,\zeta})\mathrm{d}\spatialvar,\,\,\,1\leq i \leq \dimstate,\\
(\bs{\nabla}_{\! \bparam}\objfun_{h})_{i}&=
\int_{\Omega}\left(\gamma_{1}\psi_{i}\param_{h}+\gamma_{2}\bs{\nabla}_{\! \spatialvar}\psi_{i}\cdot\bs{\nabla}_{\! \spatialvar}\param_{h}\right)\mathrm{d}\spatialvar,\,\,\,1\leq i\leq \dimparam, \\
(\bs{\nabla}^{2}_{\bstatevar,\bstatevar}\objfun_{h})_{i,j}&=
\int_{\Omega_{\text{left}}}\phi_{i}\phi_{j}\mathrm{d}\spatialvar,\,\,\,1\leq i,j \leq \dimstate, \\
(\bs{\nabla}^{2}_{\bparam,\bparam}\objfun_{h})_{i,j}&=
\int_{\Omega}\left(\gamma_{1}\psi_{i}\psi_{j}+\gamma_{2}\bs{\nabla}_{\! \spatialvar}\psi_{i}\cdot\bs{\nabla}_{\! \spatialvar}\psi_{j}\right)\mathrm{d}\spatialvar,\,\,\,1\leq i,j\leq \dimparam.
\end{align*}

\bibliographystyle{vancouver}
\bibliography{main}%

\end{document}